\documentclass[11pt, reqno]{amsart}
\usepackage[hmargin=35mm,vmargin=35mm]{geometry}
\usepackage{graphicx}
\usepackage{amssymb}

\usepackage{amsmath, amsthm, amssymb, mathabx, verbatim, setspace, enumerate,mathtools}
\usepackage[mathscr]{euscript}
\usepackage[all,cmtip]{xy}
\usepackage{etoolbox}

\newtheorem{thm}{Theorem}[section] 
\newtheorem{lem}[thm]{Lemma}  
 
\newtheorem{prop}[thm]{Proposition}  
\newtheorem{defn}[thm]{Definition}  
  
\newtheorem*{theorem*}{Theorem}

\renewcommand{\>}{\rangle}
\newcommand{\Z}{\mathbb{Z}}

\newcommand{\Q}{\mathbb{Q}}
\newcommand{\QQ}{\mathbb{Q}}

\newcommand{\PP}{\mathbb{P}}
\newcommand{\C}{\mathbb{C}}

\newcommand{\A}{\mathbb{A}}

\newcommand{\SSS}{\mathcal{S}}
\newcommand{\HH}{\mathcal{H}}

\newcommand{\AAA}{\mathbb{A}}
\newcommand{\pp}{\mathfrak{p}}

\newcommand{\p}{\mathfrak{p}}

\newcommand{\g}{\mathfrak{g}}

\newcommand{\ra}{\rightarrow}

\DeclareMathOperator{\Span}{Span}

\DeclareMathOperator{\GL}{GL}

\DeclareMathOperator{\SO}{SO}
\DeclareMathOperator{\Sp}{Sp}

\DeclareMathOperator{\Stab}{Stab}

\DeclareMathOperator{\Res}{Res}
\DeclareMathOperator{\GSpin}{GSpin}

\DeclareMathOperator{\Lie}{Lie}
\DeclareMathOperator{\CH}{CH}
\DeclareMathOperator{\Pic}{Pic}

\DeclareMathOperator{\Sym}{Sym}

\DeclareMathOperator{\SC}{SC}

\newcommand{\CC}{\mathbb{C}}
\newcommand{\ds}{\displaystyle}
\newcommand{\RR}{\mathbb{R}}
\DeclareMathOperator{\tr}{tr}
\DeclareMathOperator{\Ei}{Ei}

\title{Generating series of a new class of orthogonal Shimura varieties}
\author{Eugenia Rosu, Dylan Yott}
\address{\newline Eugenia Rosu, rosu@math.arizona.edu \newline Max Planck Institute for Mathematics, Bonn, Germany \newline Department of Mathematics, University of Arizona, Tucson, AZ 85721, USA \newline \newline Dylan Yott, dyott@math.berkeley.edu \newline Department of Mathematics, UC Berkeley, Berkeley, CA 94709, USA}
\date{}                                     
\begin{document}
\maketitle

\begin{abstract}
For a new class of Shimura varieties of orthogonal type over a totally real number field, we construct special cycles and show the the modularity of Kudla's generating series in the cohomology group. \end{abstract}

\section{\bf Introduction}

For Hilbert modular surfaces, Hirzebruch and Zagier showed in \cite{HZ} that certain generating series that have as coefficients the Hirzebruch-Zagier divisors are modular forms of weight $1$. Further inspired by this work, Gross, Kohnen and Zagier showed in \cite{GKZ} that a generating series that has Heegner divisors as coefficients is modular of weight $3/2$. This approach is unified by Borcherds in \cite{Bo}, who showed more generally the modularity of generating series with Heegner divisor classes as coefficients in the Picard group over $\QQ$.

Kudla and Millson extended the results to Shimura varieties of orthogonal type over a totally real number field and showed the modularity in the cohomology group in \cite{Ku1}, based on work from \cite{KM1}, \cite{KM2}, \cite{KM3}. This is further extended by Yuan, Zhang and Zhang in \cite{YZZ1}, who showed the modularity of the generating series in the Chow group. 

In the current paper, inspired by the above work of Kudla and Millson, we construct special cycles on a different Shimura variety of orthogonal type over a totally real number field $F$ and show the modularity of Kudla's generating series in the cohomology group. 

We consider the Shimura variety corresponding to the reductive group $\Res_{F/\QQ} G$, where $G=\GSpin(V)$ is the GSpin group for $V$ a quadratic space over a totally real number field $F$, $[F:\QQ]=d$. We choose $V$ of signature $(n,2)$ at $e$ real places and signature $(n+2,0)$ at the remaining $d-e$ places. Kudla, Millson and Yuan, Zhang, Zhang have treated the case of $e=1$, while we allow $e \in \{1, \dots, d\}$. 

If $e>1$, there is no simpler divisor case, which makes the analysis much harder. In particular, there is a very technical convergence issue that does not appear in the work of Kudla and Millson.
\bigskip

We present now the setting of the paper. For $F$ be a totally real field with real embeddings $\sigma_1, \ldots \sigma_d$, let $\A=\A_F$ be the ring of adeles of $F$ and let $V$ be a quadratic space over $F$ of signature $(n,2)$ at the infinite places $\sigma_1, \dots ,\sigma_e$ and of signature $(n+2,0)$ elsewhere. Let $G$ denote the reductive group $\GSpin(V)$ over $F$. 
We define the hermitian symmetric domain $D$ corresponding to $G$ to be: 
\[
D= D_1 \times D_2 \times \ldots \times D_e,
\]
where $D_i$ is the Hermitian symmetric domain of oriented negative definite $2-$planes in $V_{\sigma_i}=V\otimes_{\sigma_i}\RR$. 

Then $(\Res_{F/\QQ}G, D)$ is a Shimura datum and for any open compact subgroup $K$ of $G(\A_f)$, this gives us the complex Shimura variety:
\[
M_K(\CC)\simeq G(F) \backslash D \times G(\A_{f}) \slash K.
\]

For $i=1, \ldots , e$ we let $L_{D_i}$ be the complex line bundle corresponding to the points of $D_i$. We also define the projections maps $p_i: D \rightarrow D_i$ and then the line bundles $p_i^*L_{D_i} \in \Pic(D)$ descend to line bundles $L_{K, i} \in \Pic(M_K) \otimes \Q$. 


Let $W$ be a totally positive subspace of $V$, meaning that $W_{\sigma_i}=W\otimes_{\sigma_i} \RR$ is a positive subspace of $V_{\sigma_i}=V\otimes_{\sigma_i} \RR$ for all places $1\leq i\leq d$. We define $V_{W}=W^{\perp}$ to be the space of vectors in $V$ that are orthogonal to $W$, $G_{W}=\GSpin(V_{W})$ and $D_{W}=D_{W, 1}\times \dots \times D_{W, e}$ the  Hermitian symmetric domain associated to $G_W$, where $D_{W, i}$ consists of the lines in $D_i$ perpendicular to $W$. We actually have the natural identifications:
\[
G_{W}=\{g \in G: gw=w, \forall w\in V_{W}\}, \  D_W = \{(\tau_1, \ldots ,\tau_e) \in D : \langle w, \tau_i \rangle=0, \  \forall w\in W, \forall 1\leq i\leq e \},
\]
where $\left<\cdot, \cdot \right>$ is the inner product corresponding to $q_i$, the quadratic form on $V_{\sigma_{i}}$, that extends to $V_{\sigma_i}(\CC)$ by $\CC$-linearity.

Then $(\Res_{F/\QQ}G_W, D_W)$ is a Shimura datum and we have a morphism $(\Res_{F/\QQ}G_W,D_W) \rightarrow (\Res_{F/\QQ}G,D)$ of Shimura data. For $K\subset G(\AAA_f)$ an open compact subgroup and $g\in G(\AAA_f)$, we can define the complex Shimura variety:
\[
M_{gKg^{-1}, W} = G_W(F) \backslash D_W \times G_W(\A_f) \slash (gKg^{-1} \cap G_W(\A_f)).
\]

Moreover, we have an injection of $M_{gKg^{-1},W}$ into $M_K$ given by:
\[
M_{gKg^{-1},W} \rightarrow M_K, \ [\tau, h] \rightarrow [\tau, hg].
\]

We define the cycle $Z(W, g)_K$ to be the image of the morphism above. Note that $Z(W,g)_K$ is represented by the subset $D_W \times G_W(\A_f)gK$ of $D \times G(\A_f)$.

 Now let $x=(x_1, \dots, x_r)\in V(F)^r$ and let $U(x):=\Span_F\{x_1, \dots, x_r\}$ be a subspace of $V$. Then we define {\it Kudla's special cycles}:

\[
Z(x, g)_K=\begin{cases}
Z(U(x), g)_K ((-1)^e c_1(L_{K,1}^{\vee}) \ldots c_1(L_{K,e}^{\vee}) )^{r-\dim U}, & \text{if~} U(x) \text{ is totally positive,}\\
 0, & \text{ otherwise.}
\end{cases}
\]

Here $c_1$ denotes the Chern class of a line bundle. We will also use the notation $G_{x}:=G_{U(x)}$, $D_{x}:=D_{U(x)}$, $V_x:=V_{U(x)}$. Note that if $x=(x_1, \dots, x_r)\in V(F)^r$, we have $G_{x}= G_{x_1}\cap \dots \cap G_{x_r}$, as well as $D_{x}=D_{x_1}\cap \dots \cap D_{x_r}$.



\bigskip
Now we will define {\it Kudla's generating function}. For any Schwartz-Bruhat functions $\phi_f \in S(V^r(\AAA_{f}))^K$ and $g'$ in $\widetilde{\Sp}_{2r}(\AAA)$,  where $\widetilde{\Sp}_{2r}(\AAA)$ is the metaplectic cover of the symplectic group $\Sp_{2r}(\AAA)$, we define the generating series:
\[
Z(g', \phi_f) = \sum_{x \in G(F) \backslash V^r} \sum_{g \in G_x(\A_f) \backslash G(\A_f) \slash K} r(g_f')\phi_f(g^{-1}x)W_{T(x)}(g_\infty ') Z(x,g)_K.
\]
Here $r$ is the Weil representation of $\widetilde{\Sp}_{2r}(\A)\times O(V^r_{\AAA})$, where  $T(x)=\frac{1}{2}(\<x_i, x_j\>)_{1\leq i, j \leq r} \in M_r(F)$ is the intersection matrix of $x$, and $W_{T(x)}$ is the standard Whittaker function for $T(x)$.
Note that when $e=1$, for $g_f=\text{Id}$ and a careful choice of $g'_{\infty}$ we recover the generating series presented by Yuan, Zhang and Zhang in \cite{YZZ1}.

The following is the main theorem of the paper:
\begin{thm} Let $\phi_f \in S(V^r(\A_{f}))^K$ be any Schwartz-Bruhat function invariant under $K$. Then the series $[Z(g', \phi_f)]$ is an automorphic form, discrete of parallel weight $1+\frac{n}{2}$ for $g' \in \widetilde{\Sp}_{2r}(\A)$ and valued in $H^{2er}(M_K,\C)$.
\end{thm}

By modularity here we mean that, for any linear function $l:H^{2er}(M_K, \CC)\ra \CC$, the generating series obtained by acting via $l$ on the cohomology classes of the special cycles
\[
l(Z(g', \phi_f)) = \sum_{x \in G(F) \backslash V^r} \sum_{g \in G_x(\A_f) \backslash G(\A_f) \slash K} r(g_f')\phi_f(g^{-1}x)W_{T(x)}(g_\infty ') l(Z(x,g)_K).
\]
 is absolutely convergent and an automorphic form with coefficients in $\CC$ in the usual sense.

The case $e=1$ was proved by Kudla and Millson in \cite{Ku1}, based on work from \cite{KM1}, \cite{KM2}, \cite{KM3}. Yuan, Zhang and Zhang proved further in \cite{YZZ1} the modularity of $Z(g', \phi_f)$ in the Chow group. One can further conjecture that for $e>1$ the series $Z(g', \phi_f)$ is an automorphic form, discrete of weight $1+\frac{n}{2}$ for $g' \in \widetilde{\Sp}_{2r}(\A)$ valued in $\CH^{er}(M_K)_\C$. This is out of reach at the moment, but one can expect to extend the methods of Borcherds (see \cite{Bo}) to show the modularity in the Chow group.

\bigskip
 We will present now the ideas of the proof. We prove the cases $e>1$ by extending the ideas of Kudla and Millson. For each cycle $Z(x, g)$ we want to construct a Green current $\eta(x, g)$ of $Z(x, g)$ in $M_K(\CC)$. Via the isomorphism $H_{dR}^{2er}(X_K, \CC)\simeq H^{BM}_{2er(n-1)}(X_K, \CC)$, where the former is deRham cohomology while the latter is Borel-Moore cohomology, we have the identification of cohomology classes:
\[
[Z(x, g)]=[\omega(\eta(x, g))],
\]

\noindent where $\omega(\eta(x, g))$ is the Chern form corresponding to the Green current $\eta(x, g)$.

Let $x\in V(F)^r$ such that $U(x):=\Span_F\{x_1, \dots, x_r\}$ is a totally positive definite $k$-subspace of $V$. Define 
\[
x':=(x'_1,\dots, x'_k)
\] 
such that $x'_1=x_{i_1}, \dots, x'_k=x_{i_k}$ with $1\leq i_1<\dots < i_k\leq r$ the smallest indices for which $U(x')=U(x)$. 

We take the currents defined by Kudla and Millson $\eta_{0}(x'_j, \tau_i)$ of  $D_{x_j, i}$ in $D_i$, where $1\leq j\leq k$, $1\leq i\leq e$. Taking further the $*$-product of the currents $\eta_0(x_j', \tau_i)$ for $1\leq i \leq e$, we get a Green current of $D_{x, i}$ in $D_i$:
\[
\eta_1(x', \tau_i)=\eta_0(x'_1, \tau_i)*\eta_0(x'_2, \tau_i)*\dots *\eta_0(x'_k, \tau_i).
\]
Taking the pullbacks via the projections $p_i:D\ra D_i$ and taking the $*$-product, we obtain a Green current of $D_{x}$ in $D$:
\[
\eta_2(x', g)=p_1^*\eta_1(x', \tau_1)*p_2^*\eta_1(x', \tau_2)*\dots *p_e^*\eta_1(x', \tau_e).
\]

 Furthermore, we average the current  $\eta_2(x', g)$ on a lattice to get
 \[
 \ds \eta_3(x', \tau; g, h)=\sum\limits_{\gamma\in G_x(F)\setminus G(F)} \eta_2(x', \gamma \tau)1_{G_x(\AAA_f)gK}(\gamma h),
 \]
which is a Green current for $G(F)(D_{x}\times G_{x}(\AAA_f)gK/K)$ in $D\times G(\AAA_f)/K$. Showing the convergence of the sum in the definition $\eta_3(x', \tau; g, h)$ represents the most technical part of the proof and it is treated in Section \ref{conv_eta}.

As $\eta_3(x', \tau; g, h)$ is invariant under the left action of $G(F)$, $\eta_3(x', \tau; g, h)$ descends to a Green current  $\ds \eta_4(x', \tau; g, h)$ of $Z(U(x), g)_K$ in $M_K$. Here $G(F)(\tau, h)K\in M_K$.
 
Taking the Chern forms, the $*$-product turns into wedge product and the averages, as well as the pullbacks are preserved. $\omega_0(x'_j, \tau_i)$ is the Chern form of $\eta_0(x'_j, \tau_i)$ that is defined by Kudla and Millson in \cite{Ku1}, based on work from \cite{KM1}, \cite{KM2}, \cite{KM3}. Furthermore, we have
\[
\omega_1(x', \tau_i)=\omega_0(x_1, \tau_i)\wedge \dots \wedge\omega_0(x'_k, \tau_i),
\]
\[
\omega_2(x', \tau)=p_1^*\omega_1(x', \tau_1))\wedge p_2^*\omega_1(x', \tau_2)\wedge \dots \wedge p_e^*\omega_1(x', \tau_e))
\]
are the Chern forms of $\eta_1(x', \tau_i)$ and $\eta_2(x', \tau)$ respectively, and
\[
\omega_3 (x', \tau; g, h)=\sum\limits_{\gamma\in G_x(F)\setminus G(F)} \omega_2(x', \gamma \tau)1_{G_x(\AAA_f)gK}(\gamma h),
\]
\noindent is the Chern form of the Green current $\eta_3(x', \tau; g, h)$. Finally, $\omega_3 (x', \tau; g, h)$ descends to $\omega_4(x', \tau; g, h)$ corresponding to the divisor $Z(U(x), g)_K$ in $M_K$ and is the Chern form of $\eta_4(x', \tau; g, h)$.

We defined above $\omega_2, \omega_3$ and $\omega_4$ for $x' \in V(F)^k$ with $\dim U(x')=k$. We actually can extend the definitions of $\omega_2, \omega_3$ and $\omega_4$ for $x \in V(F)^r$ when $\dim U(x)<r$ as well. For $x\in V(F)^r$, if $\dim U(x)=k$, we have the equality of cohomology classes $[Z(U(x), g)]=[\omega_4(x', \tau; g, h)]$ in $H^{2ek}(M_K, \CC)$ and we can actually show further that we also have:
 \[
 [Z(x, g)]=[\omega_4(x, \tau; g, h)]
 \]
in $H^{2er}(M_K, \CC)$.
Plugging in $[\omega_4(x, \tau; g, h)]$ for the cohomology class of $[Z(x, g)]$, we take the the pullback $p^*$ of the natural projection map $p: D\times G(\AAA_f)/K\ra M_K$ and unwind the sums. Then we get:
\begin{equation}\label{gen_s}
p^*[Z(g', \phi)] = \sum_{x\in V(F)^r} r(g'_f)\phi_f(x) W_{T(x)}(g'_{\infty}) \omega_1(x, \tau).
\end{equation}

It is enough to show that (\ref{gen_s}) is an automorphic form with values in $H^{2er}(D\times G(\AAA_f)/K, \CC)$. We show this using the properties of the Kudla-Millson form on the weight of each individual $\omega_0(x, \tau_i)$,  as we can rewrite (\ref{gen_s}) as:
\[
p^*[Z(g', \phi)]=\sum_{x\in V(F)^r} r(g'_f)\phi_f(x) r(g'_{\infty}) (e^{-2\pi \tr T(x)} \omega_1(x, \tau)),
\]
and the RHS is a theta function of weight $(n+2)/2$ with values in $H^{2er}(D\times G(\AAA_f)/K, \CC)$, thus it is automorphic. 


\bigskip
{\bf Acknowledgements.} The authors would like to thank Xinyi Yuan for suggesting the problem and for very helpful discussions and insights regarding the problem. We would also like to thank the anonymous reviewer for detailed feedback. ER would also like to thank Max Planck Institute in Bonn for their hospitality, as well as to the IAS of Tsinghua University where part of the paper was written.


\section{\bf Background}

\subsection{Complex Geometry}

We will recall now some background from complex geometry (see for example \cite{CG}, \cite{GH}).

Let $X$ be a connected compact complex manifold of dimension $m$.  Suppose $Y$ is a closed compact complex submanifold of codimension $d$. Then $Y$ has no boundary and is thus a $2(m-d)$ chain in $X$. We can take the class of $Y$ to be $[Y] \in H_{2(m-d)} (X, \C).$
Note that we have the perfect pairing: 
\[
H_{2(m-d)}(X,\C) \times H_{dR}^{2(m-d)} (X,\C) \rightarrow \C,
\]
given by $\ds (Y, \eta) \rightarrow \int_Y \eta$. Thus $H_{2(m-d)}(X,\C) \simeq H_{dR}^{2(m-d)} (X,\C)^\vee$. We also have the perfect pairing: 
\[
H_{dR}^{2(m-d)} (X ,\C) \times H_{dR}^{2d} (X ,\C) \rightarrow \C,
\]
given by $\ds (\eta, \omega) \rightarrow \int_X \eta \wedge \omega$. Thus $H_{dR}^{2(m-d)} (X, \C)^\vee \simeq H_{dR}^{2d} (X, \C)$. We can compose these isomorphisms to get:
\begin{equation}\label{coh_isom1}
H_{2(m-d)} (X, \C) \simeq H_{dR}^{2d} (X, \C).
\end{equation}

For $X$ non-compact, we similarly can take the isomorphism:
\begin{equation}\label{coh_isom}
H^{BM}_{2(m-d)}(X,\C)\simeq H_{dR, c}^{2(m-d)}(X, \C)^{\vee} \simeq H_{dR}^{2d} (X, \C),
\end{equation}
where the first group is the Borel-Moore homology, which allows infinite linear combinations of simplexes, while the second group is the deRham cohomology with compact support, which uses closed differential forms with compact support.

Now for $Y$ a closed submanifold of $X$, in light of the above isomorphisms, a closed $2d$-form $\omega$ on $X$ in  $H_{dR}^{2d} (X, \C)$ represents the class $[Y]$ in $H_{2(m-d)}(X,\C)$ (respectively $H^{BM}_{2(m-d)}(X,\C)$ when $X$ non-compact), if and only if
\[
\ds \int_Y \eta = \int_X \omega \wedge \eta
\] 
for any closed $2(m-d)$ form $\eta$ on $X$.

If $X$ is not connected, we restrict the above to each of the connected components.

\subsection{Green currents and Chern forms}\label{green_back}
We recall some background on Green currents, following mainly \cite{GS}. 

Let $X$ be a quasi-projective complex manifold of dimension $m$. We define $A^{p,q}(X)$ and $A_c^{p,q}(X)$ to be the spaces of $(p,q)$-differential forms, and, respectively, $(p,q)$-differential forms with compact support. Let $D_{p,q}(X) = A_c^{p,q}(X)^*$ be the space of functionals that are continuous in the sense of deRham \cite{DR}. That is, for a sequence $\{\omega_r\} \in A^{p,q}(X)$ with support contained in a compact set $K \subset X$ and for $T \in D_{p,q}(X)$, we must have $T(\omega_r) \rightarrow 0$ if $\omega_r \rightarrow 0$, meaning that the coefficients of $\omega_r$ and finitely many of their derivatives tend uniformly to $0$. 

We also recall the differential operators:
\[
d=\partial + \overline{\partial}, \ \ d^c = \frac{1}{4 \pi i}(\partial - \overline{\partial}), \ \ dd^c = \frac{i}{2 \pi} \partial \overline{\partial}.
\]

\subsubsection{\bf Currents.} We define $D^{p,q} := D_{m-p,m-q}$ the space of $(p, q)$-currents. Then we have an inclusion $A^{p,q}(X) \rightarrow D^{p,q}(X)$ given by $\omega \rightarrow [\omega]$, where we define the current 
\begin{equation}\label{current}
[\omega](\alpha) = \int_X \omega \wedge \alpha,
\end{equation}
for any $\alpha \in A_c^{m-p,m-q}(X)$. 

For $Y\subset X$ a closed complex submanifold of dimension $p$, let $\iota:Y \hookrightarrow X $ be the natural inclusion and we also define a current $\delta_Y \in D^{p, p}(X)$ by:
\[
\delta_Y(\alpha) = \int_Y \iota^*\alpha,
\]
for any $\alpha \in A_c^{m-p,m-p}$.
\begin{defn}
A {\em Green current} for a codimension $p$ analytic subvariety $Y \subset X$ is a current $g \in D^{p-1, p-1}(X)$ such that 
\begin{equation}\label{green_current}
dd^c g + \delta_Y = [\omega_Y]
\end{equation}
for some smooth form $\omega_Y \in A^{p,p}(X)$. 
\end{defn}

This means for $\eta\in A_c^{m-p, m-p}$, we have:
\[
\int_X g dd^c \eta = \int_X \omega_Y \wedge \eta - \int_Y \eta.
\]

It implies that for a closed form with compact support $\eta$ the LHS equals $0$, and thus $\ds \int_X \omega_Y \wedge \eta = \int_Y \eta.$ Thus for $g$ a Green current of $Y$ in X, we have as cohomology classes in the isomorphism (\ref{coh_isom}):
\[
[Y]=[\omega_Y].
\]

\subsubsection{\bf Green functions and Green forms.} One natural way to obtain Green currents is from Green functions. For $Y \subset X$ a closed compact submanifold of codimension $1$, a  {\it Green function of $Y$} is a smooth function
\[
g: X \backslash Y \rightarrow \RR
\]
which has a logarithmic singularity along $Y$. This means that for any pair $(U,f_U)$ with $U \subset X$ open and $f_U: U \rightarrow \C$ a holomorphic function such that $Y \cap U$ is defined by $f_U =0$, then the function 
\[
g+\log | f_U |^2: U \backslash (Y \cap U) \rightarrow \RR
\]
extends uniquely to a smooth function on $U$. 

This definition can be extended for $Y \subset X$ a closed complex submanifold of codimension $p$ of $X$. We can define smooth forms $g_Y\in A^{p-1, p-1}(X)$ of logarithmic type along $Y$ such that  the current $[g_Y]\in D^{p-1, p-1}$ given as in (\ref{current}) by:
\[
[g_Y](\eta)= \int_{X} \eta \wedge g_Y,
\]
is a Green current. We call such smooth forms {\it Green forms of $Y$ in $X$}. We will occasionally abuse notation and use $g_Y$ for both the Green form and the Green current corresponding to $g_Y$.

\bigskip
\subsubsection{\bf Chern forms.} Now let $g$ be a Green function of $Y \subset X$, for $Y$ a divisor on $X$. For $U\subset X$ let $f_U=0$ be the local defining equation of $U\cap Y$. We define locally: 
\[
\omega_U = dd^c(g+\log|f_U|^2)
\]
By gluing together all $\omega_U$ we get a differentiable form $\omega_Y$ over $X$. We call this the {\it Chern form} associated to the Green function $g$. In general for $Y$ of codimension $p$ in $X$, for a Green form $g_Y$ of $Y$ in $X$ we call $\omega_Y$ the Chern form corresponding to $g_Y$.  





\subsubsection{\bf Star product.} Another natural way to get Green currents is by taking their $*$-product. For $Y, Z$ closed irreducible subvarieties of a smooth variety $X$ such that $Y$ and $Z$ intersect properly, let $g_Y, g_Z$ Green forms of $Y$ and $Z$, respectively. Then the $*$-product $[g_Y] * [g_Z]$ is defined by Gillet and Soul\'e in \cite{GS} to be:
\begin{equation}\label{star_product}
[g_Y] * [g_Z]=[g_Y]\wedge\delta_Z+[\omega_Y]\wedge g_Z,
\end{equation}
where $\ds [\omega_Y]\wedge g_Z(\eta)=\int_{X} \eta\wedge \omega_Y\wedge g_Z$ and $[g_Y]\wedge\delta_Z=\pi_*[\pi^*g_Y]$, where $\pi:Z \ra X$ is the embedding map. For the definition of pushforwards of currents see \cite{GS}. We can also define similarly the $*$-product $[g_Y] * G_Z$ for $g_Y$ a Green form of $Y$ and $G_Z$ a Green current for $Z$ (see \cite{GS}).
  
Moreover, from \cite{SABK} (Theorem $4$, page $50$), when $Y$ and  $Z$ have the Serre intersection multiplicity $1$, then $[g_Y]*[g_Z]$ is a Green current for $Y\cap Z$ and we have:
\begin{equation}\label{thm_S}
dd^c ([g_Y] * [g_Z]) = [\omega_Y \wedge \omega_Z]- \delta_{Y\cap Z}.
\end{equation}


\subsubsection{\bf Pullback.} Also from \cite{SABK} ($3.2$, page $50$) for $Z$ an irreducible smooth projective complex variety such that $f:Z\ra X$ is a map with $f^{-1}(Y)\neq Z$, then if $g_Y$ is a Green form of logarithmic type along $Y$, $f^*g_Y$ is a Green form of logarithmic type along $f^{-1}(Y)$. We define the pullback of currents $f^*[g_Y]:=[f^*g_Y]$ and, when the components of $f^{-1}(Y)$ have Serre intersection multiplicity $1$, the current $f^*[g_Y]$ satisfies:
\begin{equation}\label{pullback}
dd^c f^*[g_Y]+\delta_{f^{-1}(Y)}=[f^*\omega_Y].
\end{equation}






\section{\bf Construction of Green currents and Chern forms}
In this section we will construct a Green current of $Z(U, g)_K$ in $M_K$ for $U$ a totally positive subspace of $V(F)$. 

\subsection{The Shimura Variety}\label{shimura_var} Recall $\sigma_1, \dots, \sigma_d$ are the embeddings of $F$  into $\RR$ and let $(V, q)$ be a quadratic space such that $V_{\sigma_i}=V\otimes_{\sigma_i} \RR$, has signature $(n, 2)$ for $1\leq i\leq e$ and signature $(n+2, 0)$ otherwise. $V$ has the inner product given by $\left<x, y\right>=q(x+y)-q(x)-q(y)$. This can be naturally extended to $V_{\sigma_i}$ at each place $\sigma_i$ for $1\leq i\leq d$, and we denote by $q_i$ the quadratic form corresponding to this inner product.

We defined in the introduction the Hermitian symmetric domain 
\[
D=D_1\times\dots \times D_e,
\]
where $D_i$ consists of all the oriented negative definite planes in $V_{\sigma_i}$.  We can actually write explicitly the definition of $D_i$ as:
\[
D_i=\left\{v\in V_{\sigma_i}(\CC):  \left<v, v\right>=0, \, \left<v, \bar{v}\right><0 \right\}/\CC^\times \subset \mathbb{P}(V_{\sigma_i}(\CC)),
\]
where $\left<\cdot, \cdot \right>$ is the inner product corresponding to $q_i$ that extends to $V_{\sigma_i}(\CC)$ by $\CC$-linearity, and $v\mapsto \bar v$ is the involution on $V_{\sigma_i}(\CC)=V_{\sigma_i}\otimes _\RR\CC$ induced by complex conjugation on $\CC$.

We now recall the definition of $\GSpin(V)$. Let $(V, q)$ be a quadratic space over $F$ and $C(V,q) = \left( \oplus_k V^{\otimes k} \right) \slash I
$ be the Clifford algebra of $(V,q)$, where we are taking the quotient by the ideal $I = \{q(v) - v \otimes v | \  v \in V\}.$

Then $C(V,q)$ has dimension $2^{\dim(V)}$ and we have a $\Z$-grading on $T(V) =  \bigoplus_{k} V^{\otimes k} $. The map $V\ra V$, $v\ra -v$ naturally extends to an algebra automorphism $\alpha: C(V, q) \ra C(V, q)$. Then there is a natural $\Z \slash 2\Z$-grading on $C(V,q)$ given by $C(V,q) = C_0(V,q) \oplus C_1(V,q)$, where 
\[
C_i(V, q)=\{ x\in C(V, q): \alpha(x)=(-1)^i x \}, \ i=0, 1.
\]

We naturally have $V \subset C_1(V,q)$. Then we can define the $\GSpin$ group of $V$:
\[
\GSpin(V)  = \{g \in C_0(V,q)^\times | \ gVg^{-1}=V \},
\]

We denote by $G=\GSpin(V)$ and note that $G$ acts on $V$ by conjugation. The group $\Res_{F \slash \Q}G$ is reductive over $\QQ$ and the pair $(\Res_{F \slash \Q}G, D)$ is a Shimura datum. For $K\subset G(\AAA_F)$ an open compact subgroup, this gives us the complex Shimura variety:
\[
M_K(\CC)\simeq G(F)\backslash D\times G(\AAA_{f})/K.
\]

For more details on the Shimura variety $M_K$ see \cite{Sh}. 

\vspace{1 cm}

We also define the complex line bundle $L_{D_i}$ to be the restriction to $D_i$ of the tautological complex line bundle on $\mathbb{P}(V_{\sigma_i}(\CC))$. Then for the projection maps $p_i:D\ra D_i$, we get the line bundles $p_i^*L_{D_i}\in \Pic(D)$, which further descend to the line bundles $L_{K, i}\in \Pic(M_K)\otimes \QQ$ over $M_K$, defined to be:
\[
L_{K, i} = G(F) \setminus (p_i^*L_{D_i}\times G(\AAA_f)/K).
\]

\subsection{Green functions of $D_{x, i}$ in $D_i$}
We first recall how to construct a Green function of $D_{x, i}$ in $D_i$, where
\[
D_{x, i}=\{\tau_i \in D, \left<\tau, x\right>=0\}.
\]

 Let $\tau \in D_i$. It corresponds to a negative definite $2-$plane $W$ in $V_{\sigma_i}$ and we can write any $x \in V_{\sigma_i}$ as $x=x_\tau + x_{\tau^\perp}$ where $x_\tau \in W$ and $x_{\tau^\perp} \in W^\perp$.
We define:
\[
R(x,\tau) = -q_i(x_\tau), \ q_\tau(x) = q_i(x) + 2R(x,\tau).
\]
Note that this implies  $R(x, \tau) = 0$ if and only if $\tau \in D_{x, i}$. For $x\neq 0$ and $q_i(x)<0$, then $D_{x, i}$ is empty, and the statement that $R(x, \tau) = 0$ if and only if $\tau \in D_{x, i}$ is void, thus still true.

 In terms of an orthogonal basis we can write $\tau=\alpha +  \beta \sqrt{-1}$ with $\alpha, \beta \in V_{\sigma_i}$ such that $\langle \alpha, \beta \rangle=0$ and $\langle \alpha, \alpha \rangle = \langle \beta, \beta \rangle < 0$. Then $\tau$ corresponds to the negative oriented plane $W_{\tau}=\RR\alpha+\RR\beta \subset V(\RR)$, and we have:
\[
R(x, \tau) = -\frac{\langle x, \alpha \rangle^2}{\langle \alpha, \alpha \rangle} -\frac{\langle x, \beta \rangle^2}{\langle \beta, \beta \rangle}.
\]

Another important property that we use is $R(gx,g\tau)= R(x,\tau)$. This is easily seen in the definition above as the inner product is invariant under the action of $g$. 

Moreover, we show below that $-\log(R(x,\tau))$ is a Green function for $D_{x, i}$ in $D_i$:
\begin{lem}\label{green_log}
For fixed $x\in V$, $x \neq 0$, and $\tau \in D_i \backslash D_{x, i}$, the function $-\log(R(x,\tau))$ is a Green function for $D_{x, i}$ in $D_i$. 
\end{lem}
{\bf Proof:} Recall the line bundle $L_{D_i}$ is the restriction to $D_i$ of the tautological complex line bundle on $\mathbb{P}(V_{\sigma_i}(\CC))$. It has the fiber $L_\tau=\tau\CC\subset V_{\sigma_i}(\CC)$ and we have a map:
\[
s_x(\tau): L_\tau \rightarrow \C, \ v \mapsto \langle x, v \rangle.
\]
This defines an element $s_x(\tau) \in L_\tau^\vee$. As $\tau$ varies, we get a map
$$
s_x: D_i \rightarrow L_{D_i}^\vee, \ \tau \mapsto s_x(\tau).
$$
Then $s_x$ is a holomorphic section of the line bundle $L_{D_i}^\vee$. This section has a hermitian metric
$$
\left\Vert s_x(\tau) \right\Vert^2 = \frac{\lvert \langle x , v \rangle \rvert^2}{\lvert \langle v , \overline{v} \rangle \rvert},
$$
where $v \in L_\tau$ is any nonzero vector. In terms of an orthogonal basis we can write $v=\alpha +  \beta \sqrt{-1}$ such that $\langle \alpha, \beta \rangle=0$ and $\langle \alpha, \alpha \rangle = \langle \beta, \beta \rangle < 0 $. Then 
$$
\left\Vert s_x(\tau) \right\Vert^2 = \frac{\langle x, \alpha \rangle^2 + \langle x, \beta \rangle^2}{\left\vert \langle \alpha, \alpha \rangle + \langle \beta, \beta \rangle \right\vert} = -\frac{\langle x, \alpha \rangle^2}{2 \langle \alpha, \alpha \rangle} -\frac{\langle x, \beta \rangle^2}{2 \langle \beta, \beta \rangle}
$$
and also $\ds x_\tau = \frac{\langle x, \alpha \rangle}{\langle \alpha, \alpha \rangle} \alpha + \frac{\langle x, \beta \rangle}{\langle \beta, \beta \rangle} \beta.$

Computing directly gives us $R(x,\tau) = 2 \left\Vert s_x(\tau) \right\Vert^2.$ It follows by a theorem of Poincar\'e-Lelong (see Theorem 2, p. 41 in \cite{SABK}) that $-\log( R(x,\tau))$ is a Green function for $D_{x, i}$ in $D_i$.

\bigskip

For $x\in V(F)$ and $\tau\in D_i$, we have the Green function defined by Kudla and Millson (see \cite{Ku1}): 
\begin{equation}\label{f_green}
\eta(x, \tau) = f(2 \pi R(x, \tau)),
\end{equation}
where $\ds f(t) = -\Ei(-t) =\int_{t}^{\infty} \frac{e^{-x}}{x} \ dx$ is the exponential integral. Note that $\ds f(t)= -\log(t) - \gamma - \int_0^t \frac{e^{-x} -1}{x} \ dx$, where $\gamma$ is the Euler-Mascheroni constant. The function $f(t)$ is smooth on $(0,\infty)$, $f(t) + \log(t)$ is infinitely differentiable on $[0, \infty)$, and $f(t)$ decays rapidly as $t \rightarrow \infty$, thus using Lemma \ref{green_log} we easily see that $\eta(x, \tau)$ is a Green function of $D_{x, i}$ in $D_i$.

Furthermore, Kudla and Millson have constructed explicitly the Chern form $\varphi_{KM}^{(1)}(x, \tau)$ of $\eta(x, \tau)$. We recall its definition and properties in the following section.

Note that we can consider $\eta(x, \tau)$ as a restriction to $D_i$ of the Green function $f(2\pi \|s_x(v)\|^2)=f(2\pi \frac{\lvert \langle x , v \rangle \rvert^2}{\lvert \langle v , \overline{v} \rangle \rvert})$ of $\PP_x(V_{\sigma_{i}(\CC)}):=\{v\in \PP_x(V_{\sigma_{i}(\CC)}):  \left<v, x\right>=0\}$ inside $\PP(V_{\sigma_i}(\CC))$. Then the theory of Section \ref{green_back}, in particular the definition of the $*$-product, hold by restricting to $D_i$.


\subsection{The Kudla-Millson form $\varphi_{_{KM}}$}\label{Kudla_Millson}

We will now recall some results from Kudla (see \cite{Ku1}), based on previous work of Kudla and Millson (see \cite{KM1}, \cite{KM2} and \cite{KM3}). Our goal is to present explicitly the construction of the form $\varphi_{_{KM}}^{(1)}$. 

For this section we will use the notation $V_{\RR}$ for a quadratic space over $\RR$ with signature $(n,2)$, $G=\GSpin(V_{\RR})$ and $D$ the space of oriented negative $2-$planes in $V_{\RR}$. We fix a point $z_0 \in D$ and let $K = \Stab(z_0)$ be its stabilizer in $\GSpin(V_{\RR})$. Then 
\[
D \simeq G \slash K \simeq \SO(n,2) \slash (\SO(n) \times \SO(2)).
\]
Let $\g_0=\Lie(G)$ be the Lie algebra of $G$ and $\mathfrak{k}_0=\Lie(K)$ be the Lie algebra of $K$. We denote the complexifications of these lie algebras by $\g$ and $\mathfrak{k}$, respectively. We also can identify the Lie subalgebra $\p_0 \subset \g_0$ given by
\[
\p_0=\big\{\left( \begin{smallmatrix}0 & B \\ B^T & 0 \end{smallmatrix}\right): B\in M_{n\times 2}(\RR)\big\}\simeq M_{n\times 2}(\RR). 
\]
Moreover, we can give $\pp_0$ a complex structure using $J=\left( \begin{smallmatrix} 0 & 1 \\ -1 &  0 \end{smallmatrix}\right) \in \GL_2(\RR)$ acting as multiplication on the right. We denote by $\p_+$ and $\p_-$ the $\pm i$ eigenspaces of $\p$. Then we have a Harish-Chandra decomposition
\[
\g = \mathfrak{k} + \p_+ + \p_-.
\]

Moreover for the space of differential forms of type $(a,b)$ on $D$ we have an isomorphism:
\[
\Omega^{a,b}(D) \simeq [C^\infty(G) \otimes \bigwedge^{a,b} (\p^*)]^K,
\]
where on the RHS we have the wedge product $\bigwedge^{a,b} (\p^*)
=\bigwedge^a\p_{+}^*\wedge \bigwedge^b\p_{-}^*$ for $\p_{+}^*, \p_{-}^{*}$ the dual spaces of $\p_{+}$ and $\p_{-}$, respectively.

Recall that  $\widetilde{\Sp}_{2m}(\RR)$ is the metaplectic cover of $\Sp_{2m}(\RR)$, and let $K'$ be the preimage under the projection map $\widetilde{\Sp}_{2m}(\RR)\ra \Sp_{2m}(\RR)$ of the compact subgroup:
\[
\{\left(\begin{smallmatrix} A & B \\ -B & A \end{smallmatrix}\right), A+iB\in U(m)\},
\]
where $U(m)$ is the unitary group. The group $K'$ has a character $\det^{1/2}$ whose square descends to the determinant character of $U(m)$.

Then Kudla and Millson constructed a Schwartz form 
\[
\varphi^{\circ, (m)}_{_{KM}}(x, \tau) \in (\SSS(V_{\RR}^m)\otimes \Omega^{m, m}(D))^G,
\] where $\SSS(V_{\RR}^m)$ is the Schwartz space over $V_{\RR}^m$, and by invariance under $G$ we mean:
\[
\varphi^{\circ, (m)}_{_{KM}}(gx, g\tau)=\varphi^{\circ, (m)}_{_{KM}}(x, \tau).
\]

We present their result below:

\begin{theorem*}
There exists an element $\varphi^{\circ, (m)}_{_{KM}}(x, \tau) \in (\SSS(V_{\RR}^m)\otimes \Omega^{m, m}(D))^G$ with the following properties:
\begin{enumerate}[(1)]
		
	\item  For $k'\in K'$ such that $\iota(k')=\left(\begin{smallmatrix} A & B \\ -B & A \end{smallmatrix}\right)$ under the natural map $\iota: \widetilde{\Sp}_{2m}(\RR) \ra \Sp_{2m}(\RR)$, then we have:
	 \[
	 r(k')\varphi^{\circ, (m)}_{_{KM}}=(\det(k'))^{\frac{n+2}{2}} \varphi^{\circ, (m)}_{_{KM}}.
	 \]
	 
	\item $d\varphi^{\circ, (m)}_{_{KM}}=0$ i.e. for any $x\in V_{\RR}^m$, the form $\varphi^{\circ, (m)}_{_{KM}}(x, \cdot)$ is a closed $(m, m)$-form on $D$ which is $G_x$-invariant.
\end{enumerate}
\end{theorem*} 

 We define below $\varphi^{\circ, (m)}_{_{KM}}$ explicitly following \cite{Ku1}. The form $\varphi_{_{KM}}^{(m), \circ}$ is denoted by $\varphi^{(m)}$ in \cite{Ku1}. First we will construct $\varphi_{_{KM}}^{\circ, (1)}$. 
 
 Note that we have an isomorphism
\[
\ds [S(V_{\RR}) \otimes \Omega^{1,1}(D)]^G \simeq [S(V_{\RR}) \otimes \bigwedge^{1,1} \p^*]^K
\]
given by evaluating at $z_0$. Recall that we identified the Lie algebra $\ds  \p_0=\big\{\left( \begin{smallmatrix}0 & B \\ B^T & 0 \end{smallmatrix}\right): B\in M_{n\times 2}(\RR)\big\}\simeq M_{n\times 2}(\RR).$ Then we have the differential forms $\omega_{i, j}\in \Omega^1(D)=\Omega^{1, 0}(D)\oplus \Omega^{0, 1}(D)$, $1\leq i \leq n$, $1\leq j \leq 2$, defined by the function $\omega_{i, j} \in \pp_{0}^*$, $\omega_{i, j}: \pp_{0}\simeq M_{n\times 2}(\RR) \ra \RR$ given by the map $u=(u_{st})_{1\leq s \leq n, 1\leq t \leq 2} \ra u_{ij}$.

We first define for $x=(x^{(1)}, \dots, x^{(n+2)}) \in V_{\RR}$ the form $\varphi^{(1)}_{_{KM}}(x)$ that is also $G$-invariant:
\begin{equation}\label{phi_KM}
\varphi_{_{KM}}^{(1)}(x)=e^{-2\pi R(x, z_{0})}\left(\sum_{i, j=1}^{n}2x^{(i)}x^{(j)}\omega_{i, 1}\wedge \omega_{j, 2} - \frac{1}{2\pi}\sum_{i=1}^n\omega_{i, 1}\wedge\omega_{i, 2}\right).
\end{equation}

We further define $\varphi_{_{KM}}^{\circ, (1)}(x)$ to be $\ds\varphi_{KM}^{\circ, (1)}(x)=e^{-2\pi q_{z_{0}}(x)}\varphi^{(1)}_{_{KM}}(x),$ and finally, for $x=(x_1, \dots, x_m)\in V^m$ we define:
\begin{equation}\label{phi_KM_r}
\varphi_{_{KM}}^{(m)}(x)=\varphi^{(1)}_{_{KM}}(x_1)\wedge \dots \wedge\varphi^{(1)}_{_{KM}}(x_m),
\end{equation}
as well as $\varphi_{_{KM}}^{\circ, (m)}(x)=e^{-2\pi \sum\limits_{i=1}^{m}q_{z_0}(x_i)}\varphi^{(m)}_{_{KM}}(x).$

\bigskip

Recall the Green function $\eta(x, \tau)=f(2\pi R(x, \tau))$, where $x\in V(F)$ and $\tau\in D_i$. It has the important property (\cite{Ku2}, Proposition $4.10$):
\begin{equation}\label{ddc_1}
dd^c[\eta(x, \cdot)]+\delta_{D_{x, i}}=[\varphi^{(1)}_{_{KM}}(x, \cdot)],
\end{equation}
\noindent where $\varphi^{(1)}_{_{KM}} \in (S(V)\otimes \Omega^{1, 1}(D_i))^K$ is the Schwartz form defined above. This implies that $\varphi^{(1)}_{_{KM}}(x, \tau)$ is the Chern form corresponding to the Green function $\eta(x, \tau)$.
Note that (\ref{ddc_1}) is mentioned in \cite{Ku3}, Theorem 4.10 for $F=\QQ$, but holds in general for $F$ with a fixed real place $\sigma_i$ for which $V_{\sigma_i}$ has signature $(n, 2)$.




\subsection{Averaging of Green currents and their Chern forms}\label{average}
Now let $x=(x_1, \dots, x_r) \in V(F)^r$ such that $U(x)=\Span_F\{x_1, \dots, x_r\}$ is a totally positive $k$-subspace of $V(F)$, $k\leq r.$ Our goal is to construct a Green current of $Z(U(x), g)$ in $M_K$ and its corresponding Chern form.

We  define $x'=(x'_1,\dots, x'_k)$ such that $x'_1=x_{i_1}, \dots, x'_k=x_{i_k}$ and $U(x')=U(x)$. To make this uniquely defined, we pick the smallest indices $(i_1, \dots, i_k)$ for which this happens. Note further that as $U(x)=U(x')$, we also have $D_x=D_{x'}$, $V_x=V_{x'}$ and $G_x=G_{x'}$. 

For $\tau_i \in D_i$ and $x'_j \in V(F)$ for $1\leq j\leq r$, $1\leq i \leq e$,  we define as in (\ref{f_green}):
\[
f_i(x'_j, \tau_i):= f(2\pi R(x'_j,\tau_i))
\] 
that is a Green function of $D_{x'_j, i}$ in $D_i$.

We can further fix $z_{0, i} \in D_i$ for $1\leq i \leq e$ and we define the Kudla-Millson forms $\varphi^{(1)}_{_{KM}}(x_j', \tau_i) \in (\SSS(V)\otimes \Omega(D_i)^{(1, 1)})^G$ for $\tau_i\in D_i$, $x_j'\in \SSS(V)$, as in Section \ref{Kudla_Millson}, that satisfy the equation:
\begin{equation}
dd^c[f_i(x'_j, \cdot)]+\delta_{D_{x'_j, i}}=[\varphi_{_{KM}}^{(1)}(x'_j, \cdot)].
\end{equation}

As $x_1', \dots, x_k'$ are linearly independent, the submanifolds $D_{x_j', i}$ intersect properly inside $D_i$ and thus we can take the $*$-product of the Green functions $f_i(x'_j, \tau_i)$ for $1\leq j\leq k$. Denote
\[
\eta_1(x', \tau_i) = f_i(x'_1,\tau_i)*\dots *f_i(x'_k,\tau_i).
\]
Then, from (\ref{thm_S}), this is a Green current for $\ds D_{x, i}=D_{x', i}=D_{U(x'), i}=\bigcap_{j=1}^{k} D_{x'_j, i}$ in $D_i$ for $1\leq i \leq e$. 

As the star product turns into wedge product when we take the Chern forms (see (\ref{thm_S})), the Chern form associated to $\eta_1(x_j, \tau)$ is going to be:
\[
\omega_1(x', \tau_i)=\varphi^{(1)}_{_{KM}}(x'_1, \tau_i)\wedge\dots \wedge\varphi^{(1)}_{_{KM}}(x'_k, \tau_i).
\]

Note that $\omega_1(x', \tau_i)=\varphi_{_{KM}}^{(k)}(x', \tau_i)$ and thus from the definition (\ref{star_product}) of the star product, $\eta_1$ satisfies the equation:
\begin{equation}\label{ddc_r}
dd^c[\eta_1(x', \cdot)]+\delta_{D_{x, i}}=[\varphi^{(k)}_{_{KM}}(x', \cdot)].
\end{equation}


Let $p_i:D\ra D_i$ be the natural projections as before. Then, from (\ref{pullback}), $p_i^*\eta_1(x, \tau_i)$ is a Green function of $p_i^*D_{x, i}$ in $D$ and the form $p_i^*\varphi^{(k)}_{_{KM, i}}(x', \tau_i)$ satisfies:
\begin{equation}\label{ddc_pullback}
dd^c[p_i^*\eta_1(x', \cdot)]+\delta_{D_{x, i}}=[p_i^*\varphi^{(k)}_{_{KM}}(x', \cdot)]
\end{equation}

 By taking the $*$-product, we define for $\tau=(\tau_1, \dots, \tau_e) \in D\setminus D_{x}$:
\[
\eta_2(x', \tau) = p_1^*\eta_1(x',\tau_1)*\dots *p_e^*\eta_1(x',\tau_e).
\]

This is a Green current of $D_{x}$ in $D$. This follows from (\ref{pullback}), as the divisors $p_i^*D_{x, i}$ have Serre's intersection multiplicity $1$ in $D$. The Chern form of $\eta_2(x', \tau)$ is going to be:
\[
\omega_2(x', \tau)=p_1^*\omega_1(x', \tau_1)\wedge\dots \wedge p_{e}^*\omega_1(x', \tau_e),
\]
and it satisfies:
\begin{equation}\label{delta_2}
dd^c[\eta_2(x', \cdot)]+\delta_{D_{x}}=[\omega_2(x', \cdot)].
\end{equation}


\bigskip

We further take for $(\tau, h)\in D\times G(\AAA_f)$ the average of Green currents:
\[
\eta_3(x',\tau; g, h) = \sum_{\gamma \in G_{x}(F)\setminus G(F)}\eta_2(x',\gamma \tau)1_{G_{x}(\A_f)gK}(\gamma h).
\]
Note that this can be rewritten as 
\[
\ds \eta_3(x',\tau; g, h) = \sum_{\gamma \in \Gamma_h}\eta_2(\gamma^{-1}x',\tau),
\]
where $\ds \Gamma_h= G_{x}(F)\setminus G(F) \cap G_{x}(\A_f)gKh^{-1}$ is a lattice in $G(F)$. It is clear from the average that $\eta_3$ has a singularity along $G(F)(D_{x}\times G_{x}(\AAA_f)gK/K)$ in $D\times G(\AAA_f)/K$. However, note that it is not obvious that this function converges. We are actually going to prove in Section \ref{conv_eta} the following proposition:

\begin{prop}\label{conv}
Let $x\in V(F)^k$ such that $U(x)$ is a totally positive $k$-subspace of $V(F)$. Then the defining sum of $\eta_3(x,\tau; g, h)$ is absolutely convergent and $\eta_3(x,\tau; g, h)$ is a Green current of $G(F)(D_x\times G_x(\AAA_f)gK/K)$ in $D\times G(\AAA_f)/K$.
\end{prop}

This implies that $ \eta_3(x',\tau; g, h)$ is a Green current of $G(F)(D_x\times G_x(\AAA_f)gK/K)$ in $D\times G(\AAA_f)/K$. To get the Chern form we apply $dd^c$ locally and glue all the local forms using again \cite{SABK}, Theorem 4, page $50$.  This is possible due to the discussion at the end of the proof of Proposition \ref{conv} in Section \ref{conv_eta}.

Then $\eta_3$ has the Chern form:
\[
\omega_3(x', g; \tau, h)=\sum_{\gamma \in \Gamma_h } \omega_2(\gamma^{-1}x', \tau),
\]
where $\Gamma_h = G_{x'}(F)\setminus G(F) \cap G_{x'}(\AAA_f)gKh^{-1}$ as before.


\bigskip
As $\eta_3$ is invariant under the action of $G(F)$, it descends to a Green current via the projection map $p: D\times G(\AAA_f)/K \ra M_K$ to:
\[
\eta_4(x', \tau; g, h),
\]
\noindent  where $(\tau, h)$ represent the class $G(F)(\tau, h)K$ in $M_K$. The Green current condition (\ref{green_current}) is also preserved under the projection map, and the singularity is given by exactly the cycle $Z(U(x), g)_K$ inside the Shimura variety $M_K$. Thus we get:
\begin{prop} For $x'$ defined as above, $\eta_4(x',\tau; g, h)$ is a Green current of $Z(U(x), g)_K$ in $M_K$.
\end{prop}

Note that $\omega_3(x', \tau; g, h)$ descends as well to the Chern form $\omega_4(x',\tau; g, h)$ of $\eta_4(x', \tau; g, h)$. Moreover, the Chern form $\omega_3(x', \tau; g, h)$ is the pullback under the projection map $p: D\times G(\AAA_f)/K \ra M_K$ of $\omega_4(x', \tau)$: 
\[
\omega_3(x', \tau; g, h)=p^*\omega_4(x', \tau; g, h).
\]


\subsection{Extending notation}\label{ext}
In the previous section we have defined the Chern forms $\omega_2, \omega_3, \omega_4$ for $x'=(x'_1, \dots, x'_k)$ with the coordinates $x'_1, \dots, x'_k$ linearly independent. We want to extend the definition to $x=(x_1, \dots, x_k)$ in $V(F)^k$ when the coordinates $x_1, \dots, x_k$ are linearly dependent over $F$. 
In order to do that, we take $\omega_1(x, \tau_i)=\varphi^{(k)}_{_{KM}}(x, \tau_i)$, $\omega_2(x, \tau)=p_1^*\omega_1(x, \tau_1)\wedge \dots \wedge p_e^*\omega_1(x, \tau_e),$ and
 \[
\omega_3(x,\tau; g, h) = \sum_{\gamma \in G_{x}(F)\setminus G(F)}\omega_2(x,\gamma \tau)1_{G_{x}(\A_f)gK}(\gamma h).
\]
 
 We will show in Section \ref{conv_eta} in Proposition \ref{conv_omega} that $\omega_3$ is well-defined.  
\bigskip

 Also note that for $U$ a totally positive $k$-dimensional subspace of $V(F)$ we can pick any $y=(y_1, \dots, y_k)$ such that $U(y)=U$ and $\eta_4(y, \tau; g, h)$ is going to be a Green current of $Z(U, g)$ in $M_K$ with its corresponding Chern form $\omega_4(y, \tau; g, h)$.

\bigskip

 We can actually extend the definition of $\eta_2, \eta_3, \omega_2, \omega_3$ for $v\in \GL_k(F_{\infty})$ when $x=(x_1, \dots, x_k) \in V(F)^k$ such that $U(x)$ is a totally positive $k$-plane inside of $V$. We define:
\[
\eta_2(vx, \tau)= p_1^*\eta_1(v_1x,\tau_1)*\dots *p_e^*\eta_1(v_ex,\tau_e),
\]
where $v_i=\sigma_i(v) \in \GL_k(\RR)$ for $1\leq i \leq e$. Note that $G_{v_ix}=G_{x}$ and $D_{v_ix, i}=D_{x, i}$ for all $1\leq i\leq e$ and $\eta_2(vx, \tau)$ is a Green form of $D_{x}$ in $D$. 

We define further:
\[
\eta_3(vx,\tau; g, h) = \sum_{\gamma \in G_{x}(F)\setminus G(F)}\eta_2(vx,\gamma \tau)1_{G_{x}(\A_f)gK}(\gamma h),
\]
where $\eta_3(vx, \tau; g, h)$ is a Green form of $G(F)(D_x\times G_x(\AAA_f)gK/K)$ in $D\times G(\AAA_f)/K$. The proof of convergence is similar to the one for $\eta_3(x,\tau; g, h)$.

The Chern forms of $\eta_2(vx, \tau)$ and $\eta_3(vx, \tau)$ are going to be, respectively:
\[
\omega_2(vx, \tau)= p_1^*\omega_1(v_1x,\tau_1)\wedge \dots \wedge p_e^*\omega_1(v_ex,\tau_e),
\]
\[
\omega_3(vx,\tau; g, h) = \sum_{\gamma \in G_{x}(F)\setminus G(F)}\omega_2(vx,\gamma \tau)1_{G_{x}(\A_f)gK}(\gamma h).
\]

The Propositions \ref{conv} and \ref{conv_omega} extend as well for $\eta_3(vx, \tau; g, h)$ and $\omega_3(vx, \tau; g, h)$, thus they are well defined. As they are invariant under the action of $G(F)$, $\eta_3$ and $\omega_3$ further descend to the Green current $\eta_4(vx, \tau; g, h)$ of $Z(U(x), g)$ in $M_K$ that has the corresponding Chern form $\omega_4(vx, \tau; g, h)$.

Moreover, we extend the notation of $\omega_2, \omega_3$ for $x=(x_1, \dots, x_k)$ with $\dim U(x)\leq k$ by taking:
\[
\omega_2(vx, \tau)= p_1^*\omega_1(v_1x,\tau_1)\wedge \dots \wedge p_e^*\omega_1(v_ex,\tau_e),
\]
\[
\omega_3(vx,\tau; g, h) = \sum_{\gamma \in G_{x}(F)\setminus G(F)}\omega_2(vx,\gamma \tau)1_{G_{x}(\A_f)gK}(\gamma h).
\]
Propositions \ref{conv_omega} extends as well, making $\omega_3$ well-defined in general.

\subsection{Chern forms for $x=0$}

Recall that we defined in Section \ref{shimura_var} the line bundles $L_{K, i}\in \Pic(M_{K, i})\otimes \QQ$. For $x=0$, we claim that we can still define $\omega_i$ for $1\leq i\leq 4$ and the same relationships hold as in Section \ref{average}. Moreover, we are going to have:
\[
Z(0, g)=\omega_4(0, \tau).
\]

We define the Chern form $\omega_1(0, \tau_i)=(-1)^{r}\varphi^{(r)}_{_{KM}}(0, \tau_i)$. Here recall $\ds \varphi^{(1)}_{_{KM}}(0, \tau_i)=-\frac{1}{2\pi}\sum_{j=1}^n\omega_{j, 1}\wedge\omega_{j, 2} (\tau_i)$ and $\varphi^{(r)}_{_{KM}}(0, \tau_i)=\bigwedge^{r} \varphi^{(1)}_{_{KM}}(0, \tau_i)$ as defined in Section \ref{Kudla_Millson}.

We actually have:

\begin{lem} $\varphi^{(1)}_{_{KM}}(0, \tau_i)=-c_1(L_{D_i}^{\vee})$, for $1\leq i\leq e.$
\end{lem}

This is Corollary 4.12 in \cite{Ku3}. Kudla considers $F=\QQ$, but the result is unchanged for a totally real number field $F$ with a fixed embedding $\sigma_i$ into $\RR$ such that $V_{\sigma_i}$ has signature $(n, 2)$.

Thus from the lemma above we have $\omega_1(0, \tau_i)=(-1)^{r} c_1(L_{D_i}^{\vee})^{r}$. Then as before we define $\omega_2(0, \tau)=p_1^*\omega_1(0, \tau_1)\wedge \dots \wedge p_e^*\omega_1(0, \tau_e)$. Note that $\omega_2(0, \tau)=(-1)^{re}p_1^*c_1(L_{D_1}^{\vee})^{r}\wedge \dots \wedge p_e^*c_1(L_{D_e}^{\vee})^{r}$. Furthermore, as $G_0=G$, when we average over $\Gamma_h=G_0(F)\setminus (G(F) \cap G_0(\AAA_f)gKh^{-1})$ we get:
\[
\omega_3(0, \tau; g, h)=\omega_2(0, \tau).
\]
Moreover, we have as before $\omega_3(0, \tau)=p^*\omega_4(0, \tau)$, and thus 
\[
\omega_4(0, \tau)=(-1)^{re}p^*p_1^*c_1(L_{D_1}^{\vee})^r\dots p^*p_e^*c_1(L_{D_e}^{\vee})^r=(-1)^{re}c_1(L_K^{\vee})^r,
\] 
where $c_1(L_K^{\vee}):=c_1(L_{K, 1}^{\vee})\dots c_1(L_{K, e}^{\vee})$. Finally, note that $\omega_4(0, \tau)$ is exactly the cycle $Z(0, g)_K$ in $M_K$.

\bigskip


\subsection{Convergence of $\eta_3(x, \tau; g, h)$ and $\omega_3(x, \tau; g, h)$}\label{conv_eta}

Now we are ready to show the convergence of $\eta_3(x, \tau; g, h)$. More precisely, we are going to prove Proposition \ref{conv}.




\bigskip
Before we continue, we mention two short lemmas that tell us about the behavior of $R(x, \tau)$ when $\tau$ varies in a compact set in $D_i$ and $x$ varies in a lattice. 
The first lemma tells us that the quadratic forms $q_{\tau}$ bound each other:

\begin{lem}
Let $K_i \subset D_i$ be a compact set. Fix $\tau_0 \in K_i$. Then there exist $c,d >0$ such that 
\[
c q_{\tau_0}(x) \leq q_\tau(x) \leq d q_{\tau_0}(x)
\]
for all $\tau \in K_i$. 
\end{lem}
{\bf Proof:} Let $\tau \in K_i$ and $x \in V$, $x \neq 0$. Consider the function  $\psi: K_i \times \{x \in V | \ q_{\tau_0}(x)=1 \} \rightarrow \RR$, $\psi(\tau, x)=q_\tau(x).$ Since $q_{\tau_0}$ is positive definite, the set of vectors of norm $1$ is a sphere and thus compact. Hence the domain is compact and thus the image is compact, and thus bounded. Since $x \neq 0$, it must also be bounded away from $0$. Thus we can find constants $c, d$ such that:
\[
c \leq q_\tau\left(\frac{x}{\sqrt{q_{\tau_0}(x)}}\right) \leq d
\]
and $c q_{\tau_0}(x) \leq q_{\tau}(x) \leq d q_{\tau_0}(x)$ as desired. 

\bigskip

The second lemma tells us how $R(x, \tau)$ increases when $x$ varies in a lattice:

\begin{lem}\label{bounded}
 For a compact set $K_0\subset D$ and a lattice $\Gamma\subset G(F)$, there are only finitely many $\gamma \in \Gamma$ such that $R(\gamma^{-1}x, \tau_i)\leq N$ for any $\tau=(\tau_1, \dots, \tau_e) \in K_0$. More precisely, if $\dim V=n+2$, we have at most $O(N^{n/2+1})$ such $\gamma\in \Gamma$.

\end{lem}
{\bf Proof:} Fix some $\tau_0 \in K_0\cap D_i$. If for $y\in\Gamma x$ we have $R(y,\tau_i) = \frac{q_{\tau_i}(y) -a}{2} < N,$ then from the previous lemma this implies that there exists $c>0$ such that  $q_{\tau_0}(y) < \frac{a+2N}{c}.$
Thus $y$ lies in a $n+2$ dimensional sphere in $V$ of radius $\sqrt{\frac{a+2N}{c}}$. The result follows. 

\bigskip

Now we want to compute the summands of:
\begin{equation}\label{eta_3}
\eta_3(x, g; \tau, h)=\sum_{\gamma \in \Gamma_h} p_1^*\eta_1(\gamma^{-1}x, \tau_1)*p_2^*\eta_1(\gamma^{-1}x, \tau_2)*\dots *p_e^*\eta_e(\gamma^{-1}x, \tau_e),
\end{equation}
where  $\Gamma_h=G_x(F)\setminus G(F)\cap G_x(\AAA_f)gKh^{-1}$. Recall $\eta_1(x, \tau_i)=\eta_0(x_1, \tau_i)*\dots *\eta_0(x_k, \tau_i)$, where $\eta_0(x, \tau_i)=f(2\pi R(x, \tau_i))$.

We compute first the general formula for the $*$-product of $N$ Green currents:

\begin{lem}  Let $f_1, \dots, f_N$ Green forms for the cycles $Y_1, \dots, Y_N$ inside $X$, chosen such that the star product $[f_1]*\dots *[f_N]$ is well-defined. Let $\varphi_1, \dots, \varphi_N$ be their corresponding Chern forms. Then we have the $*$-product of $N$-terms:
\[
[f_1]*[f_2]*\dots *[f_N]=\sum_{j=1}^{N} \varphi_{1} \wedge \dots \varphi_{j-1}\wedge [f_{j}]\wedge \delta_{Y_{j+1}} \wedge \dots \wedge\delta_{Y_N}.
\]
\end{lem}

{\bf Proof:} We denote $\delta_{i, j}=\delta_{i}\wedge \delta_{i+1}\dots \wedge \delta_j$,  $\varphi_{i, j}=\varphi_{i}\wedge \dots \wedge \varphi_j$ for $i\leq j$ and we take $\delta_{i, j}=\varphi_{i, j}=1$ for $i>j$. We show the result by induction. For $n=2$, we have $[f_1]*[f_2]=f_1\wedge \delta_2 + \varphi_1\wedge f_2$. Assume the result is true for $n$.  Then we have: 
\[
[f_2]*[f_3]*\dots* f_{n+1} =\sum_{k=2}^{n+1}  \varphi_{2, k-1} \wedge [f_{k}]\wedge \delta_{k+1, n+1}.
\]

By definition, we have 
\[
[f_{1}]*([f_2]*[f_3]*\dots *[f_{n+1}]) =  [f_{1}]\wedge \delta_{2, n+1} +\varphi_{1}\wedge([f_2]*[f_3]*\dots *[f_{n+1}]) 
\]
\[
=  [f_{1}]\wedge (\delta_{2, n+1}) +
\sum_{k=2}^{n+1} \varphi_1\wedge\varphi_{2, k-1}\wedge [f_k]\wedge \delta_{k+1, n+1} 
\]
This is exactly $\ds \sum_{k=1}^{n+1}\varphi_{1, k-1}\wedge [f_k]\wedge \delta_{k+1, n+1}$ which finishes the proof.

\bigskip

We want to apply the above lemma to each of the $*$-products summands in (\ref{eta_3}) that define $\eta_3$: 
\[
p_1^*\eta_0(\gamma^{-1}x_1, \tau_1)*\dots p_1^*\eta_0(\gamma^{-1}x_k, \tau_1)*\dots*p_e^*\eta_e(\gamma^{-1}x_1, \tau_e)*\dots *p_e^*\eta_e(\gamma^{-1}x_k, \tau_e).
\]

Denote $f_i=p_i^*\eta_0$ and $\varphi_i=p_i^*\omega_0$. Then we get the terms:
\begin{equation}\label{term}
\sum_{i=1}^{e}\sum_{j=1}^{k} \varphi_1(\gamma^{-1}x_1, \tau_1)\wedge \dots \wedge f_i(\gamma^{-1}x_j, \tau_1)\wedge \dots \wedge\delta_{p_e^*D_{x_k}},
\end{equation}
where all the terms before $f_i$ are the smooth forms $\varphi$ and all the terms following $f_i$ are the operators $\delta$.

\bigskip

{\bf Proof of Proposition \ref{conv}:} To show the convergence of $\eta_3$, we need to show that for $\mu$ a smooth form with compact support, the integral $\ds \int_X \eta_3\wedge \mu$ converges, where $X=D\times G(\AAA_f)/K$. Note that we can cover the compact support $supp(\mu)$ of $\mu$ by finitely many open sets and in each of them we can write $\mu$ in local coordinates as a linear combination of smooth functions that are bounded inside $supp(\mu)$. 
Thus it is enough to show that the form $\eta_3$ converges to a smooth form on compacts.

We are interested in averaging the terms (\ref{term}):
\[
\sum_{i=1}^{e}\sum_{j=1}^{k} \varphi_1(y_1, \tau_1)\wedge \dots \wedge f_i(y_j, \tau_1)\wedge \dots \wedge\delta_{p_e^*D_{x_k}},
\]
 for $\tau$ inside a compact set $K_0 \subset D$,  where the average is taken over $y=(y_1, \dots, y_k)\in\Gamma_h x$. For the terms containing at least one $\delta$, the terms 
\[
 \varphi_1(\gamma^{-1}x_1, \tau_1)\wedge \dots \wedge  f_i(\gamma^{-1}x_j, \tau_1)\wedge \dots \wedge\delta_{p_e^*D_{x_k}}
\]
are nonzero only for $\tau_e \in D_{\gamma^{-1}x_k, e}$. However, this implies $R(\gamma^{-1}x_k, \tau_e)=0$ and this only happens for finitely many $\gamma \in \Gamma$ when $\tau_e \in K_0$ inside a compact from Lemma \ref{bounded}.
Thus the sum:
\[
F_1(x, \tau)=\sum_{j=1}^{k}\sum_{\substack{i=1\\ (i, j\neq (e, k))}}^{e}\sum_{\gamma\in \Gamma_h} \varphi_1(\gamma^{-1}x_1, \tau_1)\wedge \dots \wedge  f_i(\gamma^{-1}x_j, \tau_1)\wedge \dots \wedge\delta_{p_e^*D_{x_k}}
\]
is finite. This leaves the last term:
\[
F_2(x, \tau)=\sum_{\gamma\in \Gamma_h}\varphi_{1}(\gamma^{-1}x_1, \tau_1) \wedge \dots \wedge \varphi_{e}(\gamma^{-1}x_{k-1}, \tau_{e})\wedge f_{e}(\gamma^{-1}x_{k}, \tau_{e}),
\]
which we treat below in Lemma \ref{average_2}. We show that the sum $F_2(x, \tau)$ converges uniformly on compacts to a smooth form. This finishes the proof of the convergence in Proposition \ref{conv}.

Note that $F_1(x, \tau)$ is a finite sum of forms, while $F_2(x, \tau)$ is the average of wedge products of smooth forms which converges to a smooth form. 

To check the Green current condition (\ref{green_current}) is met by $\eta_3(x, \tau; g, h)$, again it is enough to check the condition on compact sets. Note first that $\tau_i \in D_{y_i}$ only for finitely many $y\in \Gamma_h x$ when $\tau$ is inside a compact set $K_0$. For $\tau_i \in  D_{y_i}$ then we have a finite sum of terms $\eta_2$ that satisfy the Green current condition (\ref{green_current}): $dd^c \eta_2(y, \tau)+\delta_{D_{y, \tau}}=[\omega_2(y, \tau)]$. For all the other terms, we do not have singularities, and as $\ds \sum_{\gamma\in \Gamma_h} \eta_2(\gamma^{-1}x, \tau)$ and all its derivatives converge to a smooth form, we can just take $dd^c$ to get 
\[
\ds dd^c \sum_{\gamma\in \Gamma_h} \eta_2(\gamma^{-1}x, \tau)=\sum_{\gamma\in \Gamma_h} dd^c\eta_2(\gamma^{-1}x, \tau)=\sum_{\gamma\in \Gamma_h} \omega_2(\gamma^{-1}x, \tau),
\] 
giving us the condition (\ref{green_current}) for $\eta_3$. Moreover, note that its Chern form is:
\[
\omega_3(x, \tau; g, h)=\sum_{\gamma\in \Gamma_h} \omega_2(\gamma^{-1}x, \tau).
\]



This finishes the proof of Proposition \ref{conv}. 

\bigskip
As promised, we show the convergence of $F_2(x, \tau)$ below:

\begin{lem}\label{average_2} The average 
\[
F_2(x, \tau; g, h)=\sum_{y\in \Gamma_h x}\varphi_{1}(y_1, \tau_1)\wedge \dots \wedge\varphi_{1}(y_k, \tau_1)\wedge \dots \wedge \varphi_{e}(y_{1}, \tau_e)\wedge\dots \wedge\varphi_{e}(y_{k-1}, \tau_{e})\wedge f_e(y_k, \tau_e)
\]
converges uniformly on compacts to a smooth form.
\end{lem}

{\bf Proof:} Let $K_0$ be a compact. We are free to discard finitely many terms from our average of the star product without affecting the convergence, so we discard the terms for which $f_e(y_k, \tau_e)=0$ on $K_0$. For $y=(y^{(1), i}, \dots, y^{(n+2), i})$ coordinates determined by the point $z_{0, i}$ in $D_{y, i}$, we recall the explicit definition of $\varphi_i(y, \tau_i)=p_i^*\varphi_{_{KM}}(y, \tau_i)$ that we presented in Section \ref{Kudla_Millson}:
\[
\varphi_i(y, \tau_i)=e^{-2\pi R(y, z_{0, i})} \left(\sum_{1\leq s, t \leq n} y^{(s), i}y^{(t), i} p_i^*(\omega_{s, 1i}\wedge\omega_{t, 2i})-\frac{1}{\pi}\sum_{1\leq s \leq n}p_i^*(\omega_{s, 1i}\wedge\omega_{s, 2i})\right).
\] 

Thus, in the average, all the terms are of the form:
\[
e^{-2\pi \sum\limits_{j=1}^{k}\sum\limits_{i=1}^{e}R(y_j, z_{0, i})} e^{2\pi R(y_k, z_{0, e})}f_e(y_k,\tau_e) \bigwedge_{i=1}^{e} \bigwedge_{\substack{j=1 \\ (i, j)\neq (e, k) }}^{k} (y_j^{(s), i} y_j^{(t), i})^{f}  p_i^*\omega_{s, 1i} \wedge  p_i^*\omega_{t, 2i} (\tau_i)
\]
The forms $p_i^*\omega_{s, 1i}, p_i^*\omega_{s, 2i}$ are smooth on $K_0$ and the values of the smooth functions representing them in local coordinates are bounded inside a compact. As they are independent of $y$, the convergence of $F_2(x, \tau)$ reduces to the convergence of: 
\[
\sum_{y \in \Gamma_h x} e^{-2\pi \sum\limits_{i=1}^{e-1}\sum\limits_{j=1}^{k} R(y_j, z_{0, i})} e^{-2\pi \sum\limits_{j=1}^{k-1}R(y_j, z_{0, e})}  f_e(y_k,\tau_e) P(y) .
\]

Here $\ds P(y)=\prod_{i=1}^{e}\prod_{\substack{j=1\\(i, j)\neq (e, k)}}^{k} \sum_{1, \leq s, t \leq n} \sum_{f=0}^{1}  (y_{j}^{(s), i} y_{j}^{(t), i})^{f}$  is a polynomial of degree $2k(e-1)$. 

Similarly, for computing the derivatives of $F_2(x, z)$ we are reduced to computing averages of the wedge products 
\[
\frac{\partial}{\partial^{R_{1, 1}}\tau_1 \partial^{S_{1, 1}} \overline{\tau_1}} \varphi_{1}(y_1, \tau_1)\wedge \dots \wedge
\frac{\partial}{\partial^{R_{1, k}}\tau_1 \partial^{S_{1, k}} \overline{\tau_1}} \varphi_{1}(y_k, \tau_1)
\wedge \dots
\]
\[
 \wedge 
\frac{\partial}{\partial^{R_{e, 1}}\tau_e \partial^{S_{e, 1}} \overline{\tau_i}} \varphi_{e}(y_{1}, \tau_e)
\wedge\dots \wedge
\frac{\partial}{\partial^{R_{e, k-1}}\tau_e \partial^{S_{e, k-1}} \overline{\tau_e}} \varphi_{e}(y_{k-1}, \tau_{e})\wedge 
\frac{\partial}{\partial^{R_{e, k}}\tau_e \partial^{S_{e, k}} \overline{\tau_e}} f_e(y_k, \tau_e). 
\]

We will break the proof in two main steps below:

{\bf Step 1: } We claim that it is enough to show that the sums:
\begin{equation}\label{conv_shorter}
\sum_{y \in \Gamma_h x} \frac{\partial}{\partial^{R_{e, k}}\tau_e \partial^{S_{e, k}} \overline{\tau_e}} f_e(y_k, \tau_e)
\end{equation}
converge for any integers $R_{e, k}, S_{e, k}\geq 0$.

In order to show this, let us compute first the partial derivatives in $\tau_i$ of the terms $\varphi(y_j, \tau_i)$ with $(j, i)\neq (k, e)$. We get:
\[
\frac{\partial}{\partial^R \tau_i \partial^S \overline{\tau_i}} \varphi(y_j, \tau_i) 
=
 e^{-2\pi R(y_j, z_{0, i})}\sum (y^{(s), i}y^{(t), i})^{f}\frac{\partial}{\partial^R \tau_i \partial^S \overline{\tau_i}}  p_i^*\omega_{s, 1i}\wedge  p_i^*\omega_{t, 2i} (\tau_i), 
\]
where $f\in \{0, 1\}$ and $1\leq s, t \leq n$. Since $ p_i^*\omega_{s, 2i}\wedge p_i^*\omega_{t, 2i}$ are smooth forms on compacts, the terms $\frac{\partial}{\partial^R\tau_i \partial^S \overline{\tau_i}}  p_i^*\omega_{s, 1i}\wedge  p_i^*\omega_{t, 2i} (\tau_i)$ are smooth as well.  Then the problem reduces to showing that the coefficients:  
\[
\sum_{y \in \Gamma_h x} e^{-2\pi \sum\limits_{i=1}^{e-1}\sum\limits_{j=1}^{k} R(y_j, z_{0, i})} e^{-2\pi \sum\limits_{j=1}^{k-1}R(y_j, z_{0, e})}  f_e(y_k,\tau_e) P(y) \frac{\partial}{\partial^{R_{e, k}}\tau_e \partial^{S_{e, k}} \overline{\tau_e}} f_e(y_k, \tau_e).
\]
converge on compacts.

 We can discard finitely many terms for which we have $R(y_j, \tau_i) \leq 1$ for any pair $(i, j)$ with $1\leq i \leq e$ and $1\leq j \leq k$. Then we can bound 
 \[
 \sum_{s, t=1}^{n} \sum_{f=0}^{1} (y_j^{(s), i} y_j^{(t), i})^f \leq (q_i(x_j)+R(y_j, \tau))^{n^2}.
 \]
 
 And thus we can further bound $\ds |P(y)|\leq C \prod\limits_{i=1}^{e} \prod\limits_{\substack{j=1 \\ (i, j)\neq (e, k)}}^{k} (q_i(x_j)+R(y_j, z_{0, i}))^{n^2}$. By discarding finitely many terms from the lattice, we can bound $e^{-2\pi R(y_j, \tau_i)}R(y_j, z_{0, i})^{m} \leq 1$, for any $1\leq  m \leq n^2$ and then
 \[
 e^{-R(y_j, z_{0, i})} (q_i(x_j) + R(y_j, z_{0, i}))^{n^2} \leq (q_i(x_j)+1)^{n^2},
 \]
which is a constant. Thus we need to show that the sums:
\[
C'\sum_{y \in \Gamma_h x} \frac{\partial}{\partial^{R_{e, k}}\tau_e \partial^{S_{e, k}} \overline{\tau_e}} f_e(y_k, \tau_e)
\]
converge for any integers $R_{e, k}, S_{e, k}\geq 0$, as claimed in (\ref{conv_shorter}).

{\bf Step 2:} Now we show the convergence of (\ref{conv_shorter}), in two parts.

\begin{enumerate}
\item  First we show the case of $\sum\limits_{y \in \Gamma_h x} f_e(y_k, \tau_e)$. We have $f_e(y_k, \tau_e)\leq \frac{e^{-2\pi R(y_k, \tau)}}{R(y_k, \tau_e)}\leq e^{-2\pi R(y_k, \tau_e)}$ for $R(y_k, \tau_e)\geq 1$, which happens for all except finitely many $y_k$'s from Lemma \ref{bounded}. Furthermore, also from Lemma \ref{bounded}, since there are at most $O(z^{\frac{n+2}{2}})$ vectors $y_k$ in our sum with $z\leq R(y_k,\tau_e) \leq z+1$, we are reduced to the convergence of
\[
\sum_{z=1}^\infty e^{-2\pi z} z^{\left(\frac{n+2}{2}\right)},
\]
which converges using the integral test.

\item Now we show the convergence of (\ref{conv_shorter}) for the partial derivatives in $\tau_e$ for the term $f_e(y_k, \tau_e)$. Note first that we can compute the derivatives:
\[
\frac{\partial}{\partial \tau_e}f_e(y_k, \tau_e)= \frac{e^{-2\pi R(y_k,\tau_e)}}{2\pi R(y_k, \tau_e)}\frac{\partial}{\partial \tau_e}R(y_k, \tau_e),
\]
\[
\frac{\partial}{\partial\overline{\tau_e}}f_e(y_k, \tau_e)=\frac{e^{-2\pi R(y_k,\tau_e)}}{2\pi R(y_k, \tau_e)}\frac{\partial}{\partial \overline{\tau_e}}R(y_k, \tau_e).
\]
We get in general terms of the form:
\[
\frac{\partial}{\partial^R\tau_e\partial^S\overline{\tau_e}}f_e(y_k, \tau_e)
=
e^{-2\pi R(y_k,\tau_e)}\sum_i \frac{ e^{-c_i R(y_k, \tau_e)}}{R(y_k, \tau_e)^{d_i}}P_i(\partial_{a_i, b_i} R),
\]
where the above is a finite sum, $P_i(\partial R, y_k)$ are polynomials in $\frac{\partial}{\partial^{a_i}\tau_e\partial^{b_i}\overline{\tau_e}}R(y_k, \tau_e)$,  and the constants $c_i, d_i$ are integers that satisfy $d_i\geq 1$,and $d_i>c_i\geq 0$. This can be easily shown by induction.
\\
Excluding the terms for which $R(y_k, \tau_e)\leq 1$, note that if we fix a basis $(e_1,\dots, e_{n+2})$ for $V_{\sigma_e}$, we have:
\[
\ds \frac{\partial}{\partial^{R}\tau_e\partial^{S} \overline{\tau_e}} R(y_k, \tau_e)
=
-\sum_{j=1}^{n+2} (y_k^{(j), e})^2\frac{\partial}{\partial^{R}\tau_e\partial^{S} \overline{\tau_e}} R(e_j, \tau_e),
\] 
thus we can further bound:
\[
|\frac{\partial}{\partial^{a}\tau_e\partial^{b} \overline{\tau_e}} R(y_k, \tau_e)|\leq M_{a, b}(q_e(x_k)+R(y_k, z_{0, e})),
\]
where $M_{a, b}$ is the upper bound of the values $\frac{\partial}{\partial^{a}\tau_e\partial^{b} \overline{\tau_e}} R(e_j, \tau_e)$ for $1\leq j\leq n+2$ and $\tau_e$ in our compact.

As $d_i>c_i$, for $R(y_k, \tau_e)\geq1$, we have $\frac{ e^{-2\pi c_i R(y_k, \tau_e)}}{(2\pi R(y_k, \tau_e))^{d_i}}<1$ and using the above bound we have more generally:
\[
|\frac{\partial}{\partial^R\tau_e\partial^S\overline{\tau_e}}f_e(y_k, \tau_e)| \leq M e^{-2\pi R(y_k,\tau_e)} \widetilde{Q}(R(y_k, \tau_e)),
\]
where $\tilde{Q}$ is a polynomial in $R(y_k, z_{0, e})$. Let $D$ be the degree of $\tilde{Q}$ and let $\widetilde{Q_0}(x):= \sum |a_n| x^n$ if $Q:=\sum a_n x^n$.


Similarly as before, we have at most $O(z^{\frac{n+2}{2}})$ values $y_k$ such that $z\leq R(y_k, \tau_e)\leq z+1$ for $\tau_e$ inside a compact, and the above convergence is equivalent to the convergence of
\[
\sum_{z=1}^{\infty}e^{-2\pi z} z^{\frac{n+2}{2}} \widetilde{Q_0}(z+1),
\]
which converges by the integral test.

\end{enumerate}



\bigskip

Now we are also going to show:

\begin{prop}\label{conv_omega} For $x=(x_1, \dots, x_k) \in V(F)^k$, the form
\[
\omega_3(x,\tau; g, h) = \sum_{\gamma \in G_{x}(F)\setminus G(F)}\omega_2(x,\gamma \tau)1_{G_{x}(\A_f)gK}(\gamma h)
\]
converges.

\end{prop}

{\bf Proof:} Note that the above statement follows for $\dim U(x)=k$ from the proof of Proposition \ref{conv}. For the general case the proof is similar to that 
of Lemma \ref{average_2}. Using the notation from Lemma \ref{average_2}, we can write:
\[
\omega_3(x,\tau; g, h)=\sum_{y \in \Gamma_h x }\varphi_1(y_1, \tau_1)\wedge\dots \wedge \varphi_1(y_k, \tau_1) \wedge\dots \wedge \varphi_1(y_1, \tau_e) \wedge \dots \wedge \varphi_e(y_k, \tau_e).
\]

Using the definition of $\varphi_i(y_j, \tau_i)$:
\[
\varphi_i(y_j, \tau_i)=e^{-2\pi R(y_j, z_{0, i})} \left(\sum_{1\leq s, t \leq n} y_j^{(s), i}y_j^{(t), i} p_i^*(\omega_{s, 1i}\wedge\omega_{t, 2i})-\frac{1}{\pi}\sum_{1\leq s \leq n}p_i^*(\omega_{s, 1i}\wedge\omega_{s, 2i})\right),
\] 
the terms $p_i^*\omega_{s, 1i} \wedge p_i^*\omega_{t, 1i}$ are independent of $y$, and we are reduced to the convergence of the coefficients:
\[
\sum_{y \in \Gamma_h x } e^{-2\pi\sum\limits_{i=1}^{e}\sum\limits_{j=1}^k R(y_j, z_{0, i})} P(y),
\]

\noindent where $\ds P(y)=\prod_{i=1}^{e}\prod_{j=1}^{k} \sum_{1\leq s, t \leq n} \sum_{f=0}^{1} (y_j^{(s), i} y_j^{(t), i})^f$. As in Lemma \ref{average_2}, we can bound:
\[
\sum_{1\leq s, t \leq n} \sum_{f=0}^{1} (y_j^{(s), i} y_j^{(t), i})^f \leq (R(y_j, z_{0, i})+q_i(x_j))^{n^2}.
\]

Moreover, for $(i, j)\neq (e, k)$, by discarding finitely many terms from the lattice we have $R(y_k, \tau_e)$ large enough and we can bound $e^{-2\pi R(y_j, \tau_i)}R(y_j, z_{0, i})^{m} \leq 1$, for any $1\leq  m \leq n^2$. Thus the convergence reduces to showing that 
 \[
  \sum\limits _{y \in \Gamma_h x } e^{-2\pi R(y_k, z_{0, e})} (R(y_k, z_{0, e})+q_e(x_k))^{n^2}
 \] 
 converges, or equivalently that any of the terms:
 \[
  \sum\limits _{y \in \Gamma_h x } e^{-2\pi R(y_k, z_{0, e})} R(y_k, z_{0, e})^m,
 \]
 converge for $1\leq m\leq n^2$. Again we have at most $O(z^{\frac{n+2}{2}})$ values $y_k$ such that $z\leq R(y_k, \tau_e)\leq z+1$ for $\tau_e$ inside a compact,  thus the above reduces to the convergence of:
 \[
  \sum\limits _{y \in \Gamma_h x } e^{-2\pi z} (z+1)^m z^{\frac{n+2}{2}}, 
  \]
which converges by the integral test. This finishes our proof.

\section{\bf Modularity of $Z(g', \phi)$}

We recall now the definition of the standard Whittaker function. Recall from Section \ref{Kudla_Millson} that we defined $\widetilde{\Sp}_{2r}(\RR)$ to be the metaplectic cover of $\Sp_{2r}(\RR)$, $K'$ the preimage under the projection map $\widetilde{\Sp}_{2r}(\RR)\ra \Sp_{2r}(\RR)$ of the compact subgroup $\{\left(\begin{smallmatrix} A & B \\ -B & A \end{smallmatrix}\right), A+iB\in U(r)\}$, where $U(r)$ is the unitary group. We also defined the character $\det^{1/2}$ on $K'$ whose square descends to the determinant character of $U(r)$.

For $(V_+, q_+)$ a quadratic space over $\RR$ of signature $(n+2, 0)$, let $\varphi^{\circ}_{+}(x_+) \in S(V_+^r)$ be the standard Gaussian:
\[
\varphi^{\circ}_{+}(x_+) = e^{-\pi \tr (x, x)_+},
\]
where $\frac{1}{2}(x, x)_+=\frac{1}{2}((x_i, x_j))_{1\leq i, j\leq r}$ is the intersection matrix of $x=(x_1, \dots, x_r)\in V_+^r$ for the inner product $(\cdot, \cdot)$ given by $q_{+}$ on $V_+$.

Then for $x\in V_+^r$ and $\beta=\frac{1}{2}(x, x)_+$ with $\beta$ in $\Sym_{r}(\RR)$, the group of symmetric $r\times r$ matrices, we define the $\beta$th "holomorphic" Whittaker function:
\[
W_{\beta}(g)=r(g)\varphi^{\circ}_{+}(x),
\]
where $g \in \widetilde{\Sp}_{2r}(\RR)$ and $r$ is the Weil representation of $ \widetilde{\Sp}_{2r}(\RR)\times O(V^r)$.

Using the Iwasawa decomposition of $\widetilde{\Sp}_{2r}(\RR)$, we can write each $g$ in the form:
\[
g=\left(\begin{smallmatrix} 1 & u \\ 0 & 1 \end{smallmatrix}\right)\left(\begin{smallmatrix} v & 0 \\ 0 & (v^T)^{-1} \end{smallmatrix}\right)k',~~ v\in \GL_r(\RR)^{+}, k'\in K',
\]
and we have:
\[
W_{\beta}(g)=\det(v)^{\frac{n+2}{4}}e^{2\pi i \tr \beta \tau} \det(k')^{\frac{n+2}{2}},
\]
where $\tau= u +  (v \cdot v^T)\sqrt{-1}$ is an element of $\HH_r$, the Siegel upper half-space  of genus $r$ (see \cite{YZZ1} for a reference).

We can extend this definition for $F_{\infty}$. For $g'=(g'_{j})_{1\leq j\leq d}\in \widetilde{\Sp}_{2r}(F_{\infty})=\prod\limits_{\sigma_j: F \hookrightarrow \RR} \widetilde{\Sp}_{2r}(\RR_{\sigma_j})$, we take:
\[
W_{\beta}(g'_{\infty})=\prod_{\sigma_j: F \hookrightarrow \RR} W_{\sigma_j(\beta)}(g'_{j}).
\]

Moreover, by writing each $g'_{j}=\left(\begin{smallmatrix} 1 & u_j \\ 0 & 1 \end{smallmatrix}\right)\left(\begin{smallmatrix} v_j & 0 \\ 0 & (v_j^T)^{-1} \end{smallmatrix}\right)k_j'$ using the Iwasawa decomposition and taking $\tau_j= u_j + i (v_j \cdot v_j^T)$ as above, we get:
\[
W_{\beta}(g'_{\infty})=\prod_{\sigma_j: F \hookrightarrow \RR} \det(v_j)^{\frac{n+2}{2}}e^{2\pi i \tr \sigma_j(\beta) \tau_j} \det(k_j')^{\frac{n+2}{2}}.
\]

Recall from the Introduction that we defined $T(x)=\frac{1}{2}(\<x_i, x_j\>)_{1\leq i, j \leq r}$ to be the intersection matrix in $M_r(F)$.  Note that for $1\leq i\leq e$ the intersection matrix $T(x)$ is different from the intersection matrix $\frac{1}{2}(x, x)_{+}$ above, for which the inner product $\left(\cdot, \cdot \right)$ is positive-definite.

We extend the definition of $W_{\beta}$ to $\sigma_j(\beta)\notin \Sym_{r}(\RR)$ for some $\sigma_j$, $1\leq j \leq e$, by taking $W_{\beta}(g'_{\infty})=0$.

\bigskip

For $g'\in \widetilde{\Sp}_{2r}(\AAA)$, $\phi\in (\SSS(V^r_{\AAA}))^K$, we defined in the introduction Kudla's generating series:
\begin{equation}\label{Z_g}
Z(g', \phi) = \sum_{x\in G(F)\setminus V(F)^r} \ \sum_{g\in G_x(\AAA_f)\setminus G(\AAA_f)/K} r(g'_f)\phi_f(g^{-1}x) W_{T(x)}(g'_{\infty}) Z(x, g)_K.
\end{equation}

We will show:

\begin{thm}\label{modular} The function $Z(g', \phi)$ is an automorphic form parallel of weight $1+n/2$ for $g'\in \widetilde{\Sp}_{2r}(\AAA)$, $\phi\in \SSS(V^r_{\AAA})$ with values in $H^{2er}(M_K, \CC)$.
 
\end{thm}


Recall that in $H^{2er}(M_K, \CC)$ we have $[Z(x, g)] = [\omega_4(x', \tau; g, h) \wedge ((-1)^ec_1(L_K^{\vee}))^{r-k}]$ as cohomology classes, where $c_1(L_K^{\vee})=c_1(L_{K, 1}^{\vee})\dots c_1(L_{K, e}^{\vee})$. We are actually going to show in Section \ref{cohomology} that $[Z(x, g)]=[\omega_4(x, \tau; g, h)]$ and we will replace in the sum (\ref{Z_g}) the cohomology class of the special cycle $Z(x, g)$ with the cohomology class of $\omega_4(x, \tau; g, h)$. We are going to show first the following expansion of the pullback of $[Z(g', \phi)]$ to $D\times G(\AAA_f)/K$:

\begin{lem}\label{lem_pullback} The pullback of the cohomology class $[Z(g', \phi)]$ to $D\times G(\AAA_f)/K$ is the cohomology class:
\[
p^*[Z(g', \phi)]= \sum_{x\in V(F)^r} r(g')\phi_f(h^{-1}x) W_{T(x)}(g'_{\infty}) \omega_2(vx, \tau),
\]
where $p:D\times G(\AAA_f)/K \ra M_K$ is the natural projection map and $g_i'=\left(\begin{smallmatrix} 1 & u_i \\ 0 & 1 \end{smallmatrix}\right)\left(\begin{smallmatrix} v_i & 0 \\ 0 & (v_i^T)^{-1} \end{smallmatrix}\right)k_i'$ is the Iwasawa decomposition of $g_i'=\sigma_i(g')$ for $1\leq i\leq d$.
\end{lem}

We claim that this will imply Theorem \ref{modular}. We will first discuss the pullback of cohomology classes in Section \ref{cohomology} and we will show Lemma \ref{lem_pullback} and Theorem \ref{modular} at the end of the section.


\subsection{Cohomology classes}\label{cohomology}

First we would like to understand better how we take the pullback of the cohomology classes $[\omega_3(x, \tau;g, h)]$ to $H^{2er}(D\times G(\AAA_f)/K, \CC)$. 

Note that for $x\in V(F)^r$ with $U(x)$ a totally positive $k$-subspace of $V$, and $g\in G(\AAA_f)$, we have the equality of cohomology classes $[Z(U(x), g)]=[\omega_4(x', g)]$ in $H^{2ek}(M_K, \CC)$ and we can take the pullback $[\omega_3(x', g)]$ to $H^{2ek}(D\times G(\AAA_f)/K, \CC)$. The pullback of $(-1)^ec_1(L_K^{\vee})$ to $H^{2}(D\times G(\AAA_f)/K, \CC)$ is $\omega_3(0, \tau)$.

We are actually going to show that the pullbacks of the Kudla cycles $Z(U(x), g) c_1(L_K^{\vee})^{r-k}$ can be represented by the cohomology class of $[\omega_3(x, g)]$ in $H^{2er}(D\times G(\AAA_f)/K, \CC)$ in the lemma below:

\begin{lem}\label{coh_class}   In $H^{2er}(D\times G(\AAA_f)/K, \CC)$ we have the equality of cohomology classes:
\[
[\omega_3(x')\wedge \omega_3(0)^{(r-k)}]=[\omega_3(x)]. 
\]
\end{lem}

To show this, we first recall from \cite{Ku1}, Lemma 7.3, how the pullback acts on the Kudla-Millson form $\varphi_{_{KM}}^{(k)}$. For $1\leq i\leq e$, recall that $(V_{\sigma_i}, q_i)$ is a quadratic space of signature $(n, 2)$.

\begin{lem}\label{pullback_KM}  Let $U\subset V_{\sigma_i}$ be a positive $k$-plane. For $y\in U$, let $\varphi_{+}^{\circ}\in \SSS(U^k)$ be the standard Gaussian $\varphi_{+}^{\circ}(y)=e^{-\pi q_i(y)}$. 
Let $\iota_U:D_{U, i} \ra D_i$ be the natural injection. Then under the pullback $\iota_{U}^*: \Omega^{k}(D_i) \ra \Omega^{k}(D_{U, i})$ of differential forms, we have:
\[
\iota_{U}^* \varphi_{_{KM}}^{(k), \circ}= \varphi_{+}^{\circ} \otimes \varphi_{_{KM, V_{U}}}^{(k), \circ},
\]
where $\varphi{_{KM, V_{U}}}^{(k), \circ} \in (\SSS(U^k)\otimes \Omega^{k, k}(D_{U, i}))^{K}$ is the Kudla-Millson form for the vector space $V_{i, U}=\left<U\right>^{\perp}$ and Hermitian symmetric domain $D_{U, i}$.
\end{lem}

For $x\in V(F)^r$ such that $U(x)$ is a totally positive $k$-subspace of $V$ we defined $x'=(x_{i_1}, \dots, x_{i_k})$. Let $x''=(x_{j_1}, \dots, x_{j_{r-k}})$ consist of the remaining components of $x$.  

Just for this section, we will use the notation $\omega_i^{(m)}(x, \tau)$ for $i=2, 3$ when $x=(x_1, \dots, x_m) \in V^m$. Using the above lemma, we are first going to show:

\begin{lem} With the above notation, the pullback of $\omega^{(r-k)}_3(x'', \tau; g, h)$ to $D_U \times G_U(\AAA_f)gK/K$ via the inclusion map $\iota: D_U \times G_U(\AAA_f)gK/K  \ra D \times G(\AAA_f)/K$ equals:
\begin{equation}\label{pullback_U3}
\iota^*\omega^{(r-k)}_3(x'', \tau; g, h)
=
\iota^*\omega^{(r-k)}_3(0, \tau; g, h)
\end{equation}

\end{lem}

{\bf Proof:} From the definition of $\varphi^{(r), \circ}_{_{KM}}$ we can write:
\begin{equation}\label{x''}
\varphi^{(r), \circ}_{_{KM}}(x)=\varphi^{(k), \circ}_{_{KM}}(x')\wedge \varphi^{(r-k), \circ}_{_{KM}}(x'').
\end{equation}
Then from Lemma \ref{pullback_KM}, for $\iota_U:D_{U, i} \ra D_i$ the natural embedding, we have $i_{U}^* \varphi_{_{KM}}^{(r-k), \circ}(x'')=(\varphi_{+}^{\circ} \otimes \varphi_{_{KM, V_{U, i}}}^{(r-k), \circ})(x'')
=
\varphi_{+}^{\circ}(x'')\varphi_{_{KM, V_{U, i}}}^{(r-k)}(0)$, as $x'' \in U^{r-k}$. Note that this implies:
\begin{equation}\label{pullback_U_KM}
i_{U}^* \varphi_{_{KM}}^{(r-k)}(x'')= \varphi_{_{KM, V_{U, i}}}^{(r-k)}(0).
\end{equation}

 We first want to pullback everything to $D$, via the projection maps $p_i:D\ra D_i$. We have the maps $\iota_U: D_U \hookrightarrow D$, $p_i: D \ra D_i$. Recall that
 \[
 D_U = D_{U, 1} \times \dots \times D_{U, e},
 \]
 and we can further define the embedding $\iota_{U, i}: D_{U, i} \hookrightarrow D_i$ and the projection map $p_{U, i}:D_U \ra D_{U, i}$. It is easy to see that $\iota_{U, i} \circ p_{U, i} = p_i \circ \iota_U$ as maps from $D_{U}$ to $D_i$, thus we also have the equality of pullbacks of differentials $\Omega^{r-k}(D_i) \ra \Omega^{r-k}(D_{U})$:
 \[
p_{U, i}^* \circ\iota_{U, i}^* = \iota_U^*\circ p_i^*.
 \]

 Then we get the equality:
\[
\iota_U^*p_i^* \varphi_{_{KM}}^{(r-k)}(x'', \tau_i)=p_{U, i}^* \circ\iota_{U, i}^*\varphi^{(r-k)}_{_{KM}}(x'', \tau_i).
\]

From (\ref{pullback_U_KM}), we have the RHS equal to $\ds p_{U, i}^*\varphi_{_{KM, V_{U, i}}}^{(r-k)}(0, \tau_i)$. Applying the same steps also for  $\varphi_{_{KM}}^{(r-k)}(0)$, we get: 
 \[
 \iota_U^*p_i^* (\varphi_{_{KM}}^{(r-k)}(0, \tau_i))
 =
 p_{U, i}^* \circ\iota_{U, i}^*(\varphi_{_{KM}}^{(r-k)}(0, \tau_i))
=
p_{U, i}^*(\varphi_{_{KM, V_{U, i}}}^{(r-k)}(0, \tau_i)).
 \]
 
 Thus we have:
 \begin{equation}\label{pullback_U_i}
 \iota_U^*p_i^* \varphi_{_{KM}}^{(r-k)}(x, \tau_i) =\iota_U^*p_i^* (\varphi_{_{KM}}^{(r-k)}(0, \tau_i))
 \end{equation}

Note that we can further take the wedge product of $\iota_U^*p_i^* \varphi_{_{KM}}^{(r-k)}(x, \tau_i)$ for $1\leq i \leq e$ to get 
\[
\iota_U^*\omega^{(r-k)}_2(x'')=\iota_U^* \bigwedge_{i=1}^{e} p_i^* \varphi^{(r-k)}(x, \tau_i)
 =
 \bigwedge_{i=1}^{e}  \iota_U^*p_i^* \varphi^{(r-k)}(x, \tau_i),
 \]
 and using (\ref{pullback_U_i}) this gives us $\iota_U^*(\omega_2^{(r-k)}(0, \tau))$. Note that this implies: 
 \begin{equation}\label{pullback_U}
 \iota_U^*\omega_2^{(r-k)}(x'')=\iota_U^*(\omega_2^{(r-k)}(0, \tau))
 \end{equation}

Finally, we are interested in the pullback of $\omega^{(r-k)}_3(x'', \tau; g, h)$ to $D_U \times G_U(\AAA_f)gK/K$ via the inclusion map $\iota: D_U \times G_U(\AAA_f)gK/K  \ra D \times G(\AAA_f)/K$. We have: 
\[
\iota^*\omega^{(r-k)}_3(x'', \tau; g, h)
=
\sum_{\gamma\in G_U(F)\setminus G(F)}\iota_U^*\omega^{(r-k)}_2(x'', \gamma\tau) 1_{G_U(\AAA_f)gK}(\gamma h),
\]
and using the pullback above for the RHS we get 
\[
\sum\limits_{\gamma\in G_U(F)\setminus G(F)}\iota_U^*\omega^{(r-k)}_2(0, \gamma\tau) 1_{G_U(\AAA_f)gK}(\gamma h),
\]
which equals $\ds\iota^*\omega^{(r-k)}_3(0, \tau; g, h)$. Thus we have $\ds \iota^*\omega^{(r-k)}_3(x'', \tau; g, h)
=
\iota^*\omega^{(r-k)}_3(0, \tau; g, h),$ which is the result of the lemma.

\bigskip

Note that using (\ref{pullback_U_KM}) and (\ref{ddc_r}) one can actually show that 
\[
[\varphi^{(r)}_{_{KM}}(x)]=[\varphi^{(k)}_{_{KM}}(x') \wedge \varphi_{_{KM}}^{(r-k)}(0)]
\]
as cohomology classes in $H^{2r}(D_i, \CC)$. 

Moreover, using (\ref{pullback_U}) and (\ref{delta_2}), one can further show that 
\[
[\omega_2^{(r)}(x)]=[\omega_2^{(k)}(x') \wedge \omega_2^{(r-k)}(0)]
\]
as cohomology classes in  $H^{2r}(D, \CC)$.

The proof of Lemma \ref{coh_class} below is based on the same principle.
\bigskip

{\bf Proof of Lemma \ref{coh_class}:} To show the equality of cohomology classes, we need to show that for a closed $(l-r, l-r)$-form $\mu$ with compact support, where $l$ is the complex dimension of $D \times G(\AAA_f)/K$, we have:
\begin{equation}\label{eq_coh}
\int\limits_{D \times G(\AAA_f)/K} \mu\wedge \omega_3^{(r)}(x)
= \int\limits_{D \times G(\AAA_f)/K} \mu\wedge \omega_3^{(k)}(x'')\wedge\omega_3^{(r-k)}(0)
\end{equation}

From (\ref{green_current}), for a closed form $\mu$, as $\mu \wedge \omega_3^{(r-k)}$ is a closed $(l-k, l-k)$-form we have: 
\[
\int\limits_{D \times G(\AAA_f)/K} \mu\wedge \omega_3^{(r)}(x)
=
 \int\limits_{D_U \times G_U(\AAA_f)gK/K}  \iota^*(\mu \wedge \omega_3^{(r-k)}(x'')).
\]

From (\ref{pullback_U3}), we have  $\iota^*(\mu \wedge \omega_3^{(r-k)}(x''))= \iota^*(\mu \wedge \omega_3^{(r-k)}(0))$, thus we get above:
 \begin{equation}\label{coh_1}
 \int\limits_{D \times G(\AAA_f)/K} \mu\wedge \omega_3^{(r)}(x)
=
 \int\limits_{D_U \times G_U(\AAA_f)gK/K}  \iota^*(\mu \wedge \omega_3^{(r-k)}(0)). 
\end{equation}

Using (\ref{green_current}) for $\mu \wedge \omega_3^{(r-k)}(0)$ we get as well:

\begin{equation}\label{coh_2}
\int\limits_{D \times G(\AAA_f)/K} \mu\wedge \omega_3^{(k)}(x') \wedge \omega_3^{(r-k)}(0)
=
\int\limits_{D_U \times G_U(\AAA_f)gK/K}  \iota^*(\mu \wedge \omega_3^{(r-k)}(0)).
\end{equation}
Combining the two equations (\ref{coh_1}) and (\ref{coh_2}) we get (\ref{eq_coh}).

\bigskip

{\bf Remarks on $\omega_3(vx)$ and $\omega_4(vx)$.} We follow up with some remarks regarding $\omega_3(vx, \tau; g, h)$ and $\omega_4(vx, \tau; g, h)$ when $v\in \GL_r(F_{\infty})$ and $x\in V(F)^r$ with $U(x)$ totally positive definite $k$-subspace of $V(F)$. We have defined them in Section \ref{ext}. Lemma \ref{coh_class} extends easily for  $\omega_3(vx, \tau; g, h)$ and $\omega_4(vx, \tau; g, h)$ and we have as cohomology classes in $H^{2er}(D\times G(\AAA_f)/K, \CC)$:
\[
[\omega_3(vx, \tau; g, h)]=[\omega_3((vx)', \tau; g, h) \wedge \omega^{(r-k)}_3(0, \tau)].
\]

As actually $\omega_3((vx)')$ represents the same cohomology class as the preimages of $Z(U(vx), g)$  in $D\times G(\AAA_f)/K$, and as $Z(U(x), g)=Z(U(vx), g)$, we have:

\begin{lem}
\begin{enumerate}[(i)]
\item As cohomology classes in $H^{2er}(D\times G(\AAA_f)/K, \CC)$, we have:
\begin{equation}\label{coh_v}
[\omega_3(vx, \tau; g, h)]=[\omega_3(x, \tau; g, h)].
\end{equation}

\item Noting that (\ref{coh_v}) descends to $M_K$, we also have as cohomology  classes in $H^{2er}(M_K, \CC)$:
\begin{equation}\label{coh_v4}
[\omega_4(vx, \tau; g, h)]=[\omega_4(x, \tau; g, h)].
\end{equation}

\end{enumerate}
\end{lem}


\bigskip

{\bf Proof of modularity:} We will finish below the proofs of Lemma \ref{lem_pullback} and Theorem \ref{modular}.

{\bf Proof of Lemma \ref{lem_pullback}:} The pullback to $D\times G(\AAA_f)/K$ of $\omega_4(x', \tau)$ is $\omega_3(x', \tau)$ and $\omega_3(0, \tau)$ is the pullback of $(-1)^{er}c_1^r(L_K^{\vee})=Z(0, g)$. Then in (\ref{Z_g}) we can write:

\[
p^*[Z(g', \phi)] = \sum_{x\in G(F)\setminus V(F)^r} \ \sum_{g\in G_x(\AAA_f)\setminus G(\AAA_f)/K} r(g', g)\phi_f(x) W_{T(x)}(g'_{\infty}) [\omega_3(x', \tau; g, h) \wedge \omega_3^{(r-k)}(0)].
\]

Furthermore, from Corollary \ref{coh_class} we have $[\omega_3(x', g; \tau, h) \wedge \omega_3^{(r-k)}(0, \tau)]=[\omega_3(x, \tau; g, h)]$ as classes in $H^{2er}(D\times G(\AAA_f)/K, \CC)$. From (\ref{coh_v}) we also have the equality of cohomology classes $[\omega_3(x, \tau; g, h)]=[\omega_3(vx, \tau; g, h)].$ Thus we get:
\[
p^*[Z(g', \phi)] = \sum_{x\in G(F)\setminus V(F)^r} \ \sum_{g\in G_x(\AAA_f)\setminus G(\AAA_f)/K} r(g', g)\phi_f(x) W_{T(x)}(g'_{\infty}) [\omega_3(vx, g; \tau, h)].
\]

 By plugging in the definition $\omega_3(vx, \tau; g, h)=\sum\limits_{\gamma \in G_x(F)\setminus G(F)} \omega_2(vx, \gamma\tau) 1_{G_x(\AAA_f)gK}(\gamma h)$, we get the cohomology class $p^*[Z(g', \phi)]$ equal to the cohomology class of:
\[
 \sum_{x\in G(F)\setminus V(F)^r} \sum_{g\in G_x(\AAA_f)\setminus G(\AAA_f)/K} r(g'_f, g)\phi_f(x) W_{T(x)}(g'_{\infty}) \sum_{\gamma \in G_x(F)\setminus G(F)}\omega_2(vx, \gamma\tau) 1_{G_x(\AAA_f)gKh^{-1}}(\gamma).
\]
We will unwind the sum below to get the result of the lemma. We interchange the summations to get:
\[
 \sum_{x\in G(F)\setminus V(F)^r} \sum_{\gamma \in G_x(F)\setminus G(F)} \sum_{g\in G_x(\AAA_f)\setminus G(\AAA_f)/K} r(g_f', g)\phi_f(x) W_{T(x)}(g'_{\infty}) \omega_2(vx, \gamma\tau) 1_{G_x(\AAA_f)gK}(\gamma h).
\]
Note that $1_{G_x(\AAA_f)gK}(\gamma h)\neq 0$ iff $\gamma h \in G_x(\AAA_f)gK$, or equivalently if $g\in G_x(\AAA_f)\gamma hK$, and since we are summing for $g \in G_x(\AAA_f)\setminus G(\AAA_f)/K$, we can replace $g$ by $\gamma h$ everywhere and get:
\[
p^*[Z(g', \phi)] = \sum_{x\in G(F)\setminus V(F)^r} \sum_{\gamma \in G_x(F)\setminus G(F)} r(g'_f, \gamma h)\phi_f(x) W_{T(x)}(g'_{\infty}) \omega_2(vx, \gamma\tau).
\]
Since the action of $G(\AAA_f)$ on $\phi$ is given by $r(g_f', \gamma h)\phi_f(x)= r(g_f')\phi_f(h^{-1}\gamma^{-1} x)$ and $\omega_2(vx, \gamma\tau)=\omega_2(\gamma^{-1}vx, \tau)=\omega_2(v(\gamma^{-1}x), \tau)$, then we have:
\[
p^*[Z(g', \phi)] = \sum_{x\in V(F)^r} r(g_f')\phi_f(h^{-1}x) W_{T(x)}(g'_{\infty}) \omega_2(vx, \tau),
\]
which gives us the result of the lemma. 

\vspace{1 cm}

{\bf Proof of Theorem \ref{modular}:} We would like to rewrite the sum of Lemma \ref{lem_pullback}:
 \[
p^*[Z(g', \phi)]= \sum_{x\in V(F)^r} r(g')\phi_f(h^{-1}x) W_{T(x)}(g'_{\infty}) \omega_2(vx, \tau)
\]
and first show that this sum is automorphic with values in $H^{2er}(D\times G(\AAA_f)/K, \CC)$.

We recall the Iwasawa decomoposition of $g'=(g'_{i})_{1\leq i\leq d}\in \widetilde{\Sp}_{2r}(F_{\infty})$ to be $g'_{i}=\left(\begin{smallmatrix} 1 & u_i \\ 0 & 1 \end{smallmatrix}\right)\left(\begin{smallmatrix} v_i & 0 \\ 0 & (v_i^T)^{-1} \end{smallmatrix}\right)k_i'$, where $v_i\in \GL_r(\RR_{\sigma_i})^{+}$, $k_i'\in K_i'$.

Recall that we have, for $1\leq i \leq e$, $\omega_1(x, \tau_i)=\varphi^{(r)}_{_{KM}}(x, \tau_i)$ and $\omega_2(x, \tau)=p_1^*\omega_1(x, \tau_1)\wedge\dots \wedge p_e^*\omega_1(x, \tau_e).$ From the property (1) of the Theorem of Kudla and Millson we presented in Section \ref{Kudla_Millson}, we have
\[
r(k'_{i})\varphi^{(r),\circ}_{_{KM}}=\det(k'_i)^{\frac{n+2}{2}}\varphi^{(r), \circ}_{_{KM}},
\]
\noindent where $\varphi^{(r), \circ}_{_{KM}}(x, \tau_i)=e^{-2\pi \tr \sigma_i(T(x))}\varphi_{_{KM}}(x, \tau_i)$. Using the Weil representation this easily extends to:
\[
r(g'_i)\varphi^{(r), \circ}_{_{KM}}(x, \tau_i)
=
\det(v_i)^{\frac{n+2}{2}}\det(k'_i)^{\frac{n+2}{2}}e^{-2\pi \tr T(\sigma_i(x))(u_i+i v_i \cdot v_i^T)}\varphi^{(r)}_{_{KM}}(v_ix, \tau_i).
\]

We take the pullback to $D$ via the projection maps $p_i:D \ra D_i$. We denote $\varphi_{i}(x, \tau_i)=p_i^*\varphi^{(r)}_{_{KM}}(x, \tau_i)$ and $\varphi^{\circ}_{i}(x, \tau_i)=e^{-2\pi \tr \sigma_i(T(x))}\varphi_{i}(x, \tau_i)$ and thus we also have:  
\[
r(g_i')(\varphi_{i}^{\circ}(x, \tau_i))
=
\det(v_i)^{\frac{n+2}{2}}\det(k'_i)^{\frac{n+2}{2}}e^{-2\pi \tr T(\sigma_i(x))(u_i+i v_i \cdot v_i^T)}\varphi_{i}(v_ix, \tau_i).
\]

Note that on the RHS we got $W_{\sigma_i(T(x))}(g_i')\varphi_{i}(v_ix, \tau_i)$, thus:
\[
r(g_i')(\varphi_{i}^{\circ}(x, \tau_i))
=
W_{\sigma_i(T(x))}(g_i')\varphi_{i}(v_ix, \tau_i).
\]
 Furthermore, as we can rewrite
\[
W_{T(x)}(g'_{\infty})\varphi_{1}(v_1x, \tau_1)\wedge\dots \wedge \varphi_{e}(v_ex, \tau_e)
=
\]
\[
(W_{\sigma_1(T(x))}(g'_1)\varphi_{1}(v_1x, \tau_1)\wedge\dots \wedge W_{\sigma_e(T(x))}(g'_{e})\varphi_{e}(v_ex, \tau_e))\prod_{i=e+1}^{d}W_{\sigma_i(T(x))}(g'_i),
\]
we get:
\[
W_{T(x)}(g'_{\infty})\varphi_{1}(v_1x, \tau_1)\wedge\dots \wedge \varphi_{e}(v_ex, \tau_e)
=
r(g'_{\infty})\phi^{\circ}(x, \tau),
\]
\noindent where $\phi^{\circ}(x, \tau)=\varphi_{1}^{\circ}(x, \tau_1)\wedge\dots \wedge \varphi^{\circ}_{e}(x, \tau_e)\prod\limits_{i=e+1}^{d}\varphi_{0, i}(x)$. Recall that for $i\geq e+1$, $W_{T(\sigma_i(x))}(g_i)=r(g_i)\varphi_{0, i}(x)$. Here $\varphi_{0, i}(x)=e^{-\pi \tr T(\sigma_i(x))}$ is the standard Gaussian, as $(V_{\sigma_i}, q_i)$ is positive definite for $i\geq e+1$.

Going back to the sum of Lemma \ref{lem_pullback}, we thus get:
\[
p^*[Z(g', \phi)]= \sum_{x\in V(F)^r} r(g_f')\phi_f(h^{-1}x) r(g'_{\infty})\phi^{\circ}(x, \tau),
\]
and this is a theta function of weight $(n+2)/2$ with values in the cohomology group $H^{2er}(D\times G(\AAA_f)/K, \CC)$. This means that for any linear functional $l:H^{2er}(D\times G(\AAA_f)/K, \CC) \ra \CC$ acting on the cohomology part of $\phi^{\circ}(x, \tau)$, the generating series:
\[
l(p^*[Z(g', \phi)])=\sum_{x\in V(F)^r} r(g_f')\phi_f(h^{-1}x) r(g'_{\infty})l(\phi^{\circ}(x, \tau))
\]
is a theta function of weight $(n+2)/2$. Note that this series is obtained by unwinding: 
\[
p^*[Z(g', \phi)]
=
 \sum_{x\in G(F)\setminus V(F)^r} \ \sum_{g\in G_x(\AAA_f)\setminus G(\AAA_f)/K} r(g', g)\phi_f(x) W_{T(x)}(g'_{\infty}) l(\omega_3(x, g)).
\]

Denote 
\[
\ds Z_0(g', \phi)=\sum_{x\in G(F)\setminus V(F)^r} \ \sum_{g\in G_x(\AAA_f)\setminus G(\AAA_f)/K} r(g', g)\phi_f(x) W_{T(x)}(g'_{\infty}) \omega_3(x, g).
\]
For the the natural projection $p:D\times G(\AAA_f)/K \ra M_K$, recall the pullback $p^*:\Omega^{2er}(M_K) \ra \Omega^{2er}(D\times G(\AAA_f)/K)$, which further descends to the cohomology groups $p^*:H_{dR}^{2er}(M_K) \ra H_{dR}^{2er}(D\times G(\AAA_f)/K)$ and the map is an injection.

We denote by $\SC^{2er}(M_K)$ the subspace of $H_{dR}^{2er}(M_K)$ generated by the classes $[\omega_{4}(x, g)]$ and by $\SC^{2er}(D\times G(\AAA_f)/K)$ the subspace of $H_{dR}^{2er}(M_K)$ generated by the classes $[\omega_{3}(x, g)]$. 
Then the above pullback map restricts to $p^*:\SC^{2er}(M_K) \ra \SC^{2er}(D\times G(\AAA_f)/K)$ and it is an injection.


Then for any linear functional $l$ of $\SC^{2er}(M_K)$, we can just define the linear functional $\widetilde{l}$ on $\SC^{2er}(D\times G(\AAA_f)/K)$ given by $\widetilde{l}(p^*[\omega])=\widetilde{l}([\omega])$, and thus $\widetilde{l}(Z_0(g', \phi))=l([Z(g', \phi)])$ is automorphic. Thus $[Z(g', \phi)]$ is a theta function valued in $H^{2er}(M_K)$.

 We can also easily check the weight of the theta function by computing $r(k')\phi^{\circ}(x, \tau)
=
r(k'_1)\varphi_{1}^{\circ}(x, \tau_1)\wedge\dots \wedge r(k'_e)\varphi^{\circ}_{e}(x, \tau_e)\prod\limits_{i=e+1}^{d}r(k'_i)\phi_{0, i}(x)$ which gives us the factor $\det(k'_i)^{\frac{n+2}{2}}$ at each place $i$.

\end{document}